\documentclass[12pt]{article}
\usepackage{amscd}
\usepackage[all]{xy}
\usepackage{CJK}
\usepackage{CJKnumb}
\usepackage{graphicx}
\usepackage{amsfonts}
\usepackage{amsmath}
\usepackage{amssymb}
\usepackage{fancyhdr}
\usepackage{indentfirst}
\usepackage{titlesec}
\usepackage{mathrsfs}
\usepackage{bbm}
\usepackage{dsfont}
\usepackage{amsthm} 
\usepackage{mathrsfs}

\newcommand{\cR}{\mathcal{R}}
\newcommand{\cE}{\mathcal{E}}

\newcommand{\cL}{\mathcal{L}}

\newcommand{\cA}{\mathcal{A}}

\newcommand{\cC}{\mathcal{C}}

\newcommand{\cB}{\mathcal{B}}
\newcommand{\cK}{\mathcal{K}}
\newcommand{\cH}{\mathcal{H}}
\newcommand{\cI}{\mathcal{I}}
\newcommand{\cJ}{\mathcal{J}}
\newcommand{\cS}{\mathcal{S}}

\theoremstyle{definition}  \newtheorem{Def}{Definition}
\theoremstyle{plain}  \newtheorem{Thm}{Theorem} \theoremstyle{plain}
\newtheorem{cor}{Corollary}
 \theoremstyle{remark} \newtheorem{Rek}{Remark}
 \theoremstyle{plain}\newtheorem{Lem}{Lemma}
 \theoremstyle{plain}
\newtheorem{Prop}{Proposition}
\newtheorem{Nota}{Notation}
\theoremstyle{remark}

\begin{document}

\title{\textbf{Higher Dimensional Homology Algebra V:Injective Resolutions and Derived 2-Functors in ($\cR$-2-Mod)}}
\author{Fang HUANG, Shao-Han CHEN, Wei CHEN, Zhu-Jun ZHENG\thanks{Supported in part by NSFC with grant Number
10971071
 and Provincial Foundation of Innovative
Scholars of Henan.} }
\date{}
 \maketitle

\begin{center}
\begin{minipage}{5in}
{\bf  Abstract}: In this paper, we will construct the injective
resolution of any $\cR$-2-module, define the right derived
2-functor, and give some related properties of the derived 2-functor
in ($\cR$-2-Mod).

{\bf{Keywords}:} $\cR$-2-Module; Injective Resolution; Right Derived
2-Functor
\\
\end{minipage}
\end{center}
\maketitle \hspace{1cm}

\section{Introduction}
This is the fifth paper of the series of higher dimensional homology
algebra.
In \cite{4}, We gave the definition of $\cR$-2-modules and proved
that the 2-category ($\cR$-2-Mod) is an abelian 2-category in a
different way from M.Dupont's(\cite{2}). Based on the works of A.del
R\'{\i}o, J. Mart\'{\i}nez-Moreno and E. M. Vitale\cite{11}, we gave
the constructions of projective resolutions of symmetric 2-group and
$\cR$-2-module, defined the left derived 2-functors and gave the
fundamental properties of derived 2-functors in (2-SGp) and
($\cR$-2-Mod), respectively \cite{27,28}. In this paper, we
construct an injective resolution of any $\cR$-2-module in the
2-category ($\cR$-2-Mod) and prove that it is unique up to 2-cochain
homotopy (Proposition 2 and Theorem 1). These results make it
possible to define right derived 2-functor in ($\cR$-2-Mod).

When we proved (2-SGp) and ($\cR$-2-Mod) have enough projective
objects (\cite{14}), T.Pirashvili gave injective enough of (2-SGp)
and ($\cR$-2-Mod) in \cite{20} and he said
that "abelian 2-category $\mathfrak{SCG}$ (as well as the abelian
2-category of 2-modules over a 2-ring) has enough projective and
injective objects"(\cite{24}).
He also proved that the  2-category $\mathfrak{SCG}$ is 2-equivalent
to the 2-category of 2-modules over $\Phi$, where $\Phi$ is a
symmetric categorical group, whose objects are integers, morphism
from $n$ to $m$ is $\delta_{n,m}$, i.e. if $n\neq m$, there is no
morphism; otherwise, is identity, for $n, m\in \mathds{Z}$.
So we just consider the right derived 2-functor in 2-category
($\cR$-2-Mod).

In our coming papers, we shall define $\cE$xt 2-functor and the
spectral sequence in an abelian 2-category, try to give the relation
between $\cE$xt 2-functor and the extension of 2-modules. Meantime,
we shall discuss the representation of Lie 2-algebra and its
extension based on the theory of higher dimensional homology
algebra.

This paper is organized as follows. In section 2, we give some basic
facts on $\cR$-2-modules such as the relative (co)kernel, relative
2-exact appeared in \cite{11,2,6}. The cohomology $\cR$-2-modules of
a complex of $\cR$-2-modules appear in this section, too. We give
the definition of 2-cochain homotopy of two morphisms of complexes
in ($\cR$-2-Mod) like cochain homotopy in 1-dimensional case, and
prove that it induces equivalent morphisms between cohomology
$\cR$-2-modules. In section 3, we mainly give the definition of
injective resolution of $\cR$-2-module and its
construction(Proposition 2). In section 4, after the basic
definitions of 2-functors between abelian 2-category
($\cR$-2-Mod)(\cite{2,4}), we define the right derived 2-functor and
obtain our main result Theorem 2.

\begin{Nota}
 An $\cR$-2-module we mentioned in this paper is
$(\cA,I,\cdot,a,b,i,z)$, where $\cA$ is a symmetric 2-group with
$\cR$-2-module structure $\cdot$, $I$ is the unit object under
$\cdot$, $a,b,i,z$ are natural isomorphisms satisfying canonical
properties \cite{4}, and $\cR$ is a  2-ring(\cite{2,4,5}).
\end{Nota}

\section{Preliminary}
In this section, we review the constructions of relative kernel and
cokernel(\cite{28}), give the
definition of relative 2-exactness of a sequence and the cohomology
$\cR$-2-modules of a complex of $\cR$-2-modules in
($\cR$-2-Mod)(also see \cite{2} for general defintions), and then
give the definition of 2-cochain homotopy and its related result
similar to the 1-dimensional case \cite{16,18,7,10}.
The relative kernel
$(Ker(F,\varphi),e_{(F,\varphi)},\varepsilon_{(F,\varphi)})$ of a
sequence $(F,\varphi,G):\cA\rightarrow\cB\rightarrow\cC$ in
($\cR$-2-Mod) is a symmetric 2-group consisting of:

$\cdot$ An object is a pair $(A\in obj(\cA),a:F(A)\rightarrow 0)$
such that the following diagram commutes
\begin{center}
\scalebox{0.9}[0.85]{\includegraphics{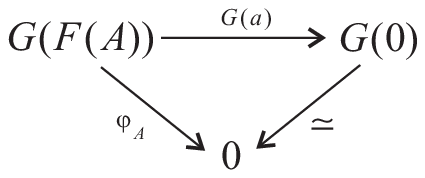}}
\end{center}

$\cdot$ A morphism $f:(A,a)\rightarrow (A^{'},a^{'})$ is a morphism
$f:A\rightarrow A^{'}$ in $\cA$ such that the following diagram
commutes
\begin{center}
\scalebox{0.9}[0.85]{\includegraphics{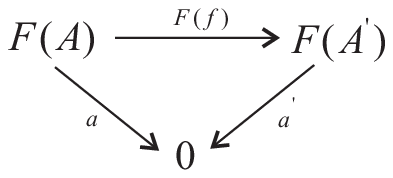}}
\end{center}

$\cdot$ $\cR$-2-module structure is induced by the $\cR$-2-module
structure of $\cA$(more details see \cite{28}).

$\cdot$ The faithful functor
$e_{(F,\varphi)}:Ker(F,\varphi)\rightarrow\cA$ is defined by
$e_{(F,\varphi)}(A,a)=A$, and the natural transformation
$\varepsilon_{(F,\varphi)}:F\circ e_{(F,\varphi)}\Rightarrow 0$ by
$(\varepsilon_{(F,\varphi)})_{(A,a)}=a$.

The relative cokernel
$(Coker(\varphi,G),p_{(\varphi,G)},\pi_{(\varphi,G)})$ of a sequence
$(F,\varphi,G):\cA\rightarrow\cB\rightarrow\cC$ in ($\cR$-2-Mod) is
a symmetric 2-group consisting of:

$\cdot$  Objects are those of $\cC$.

$\cdot$ A morphism from $X$ to $Y$ is an equivalent class of the
pair $(B,f):X\rightarrow Y$ with $B\in obj(\cB)$ and $f:X\rightarrow
G(B)+Y$. Two morphisms $(B,f),(B^{'},f^{'}):X\rightarrow Y$ are
equivalent if there is $A\in obj(\cA)$ and $a:B\rightarrow
F(A)+B^{'}$ such that the following diagram commutes
\begin{center}
\scalebox{0.9}[0.85]{\includegraphics{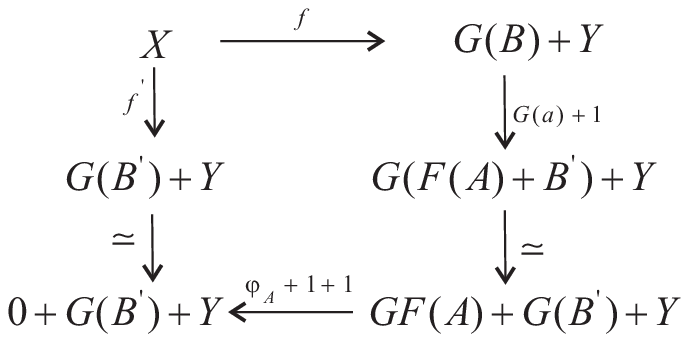}}
\end{center}

$\cdot$ $\cR$-2-module structure is induced by the $\cR$-2-module
structure of $\cC$(more details see \cite{28}).

$\cdot$ The essentially surjective functor
$p_{(\varphi,G)}:\cC\rightarrow Coker(\varphi,G)$ is defined by
$p_{(\varphi,G)}(X)=X$, and the natural transformation
$\pi_{(\varphi,G)}: p_{(\varphi,G)}\circ G\Rightarrow 0$ by
$(\pi_{(\varphi,G)})_{B}=1_{G(B)}$.



\begin{Def}(\cite{2})
Consider the following diagram in ($\cR$-2-Mod)
\begin{center}
\scalebox{0.9}[0.85]{\includegraphics{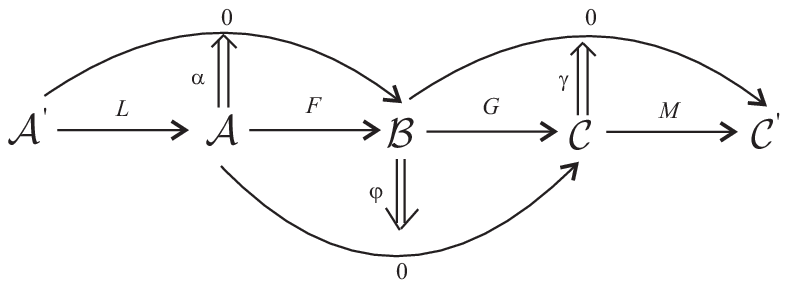}}
\end{center}
with $\alpha$ compatible with $\varphi$ and $\varphi$ compatible
with $\gamma$. By the universal property of the relative cokernel
$Coker(\alpha,F)$, we get a factorization $(G^{'},\varphi^{'})$ of
$(G,\varphi)$ through $(p_{(\alpha,F)},\pi_{(\alpha,F)})$. By the
cancellation property of $p_{(\alpha,F)}$, we have a 2-morphism
$\overline{\gamma}$ as in the following diagram
\begin{center}
\scalebox{0.9}[0.85]{\includegraphics{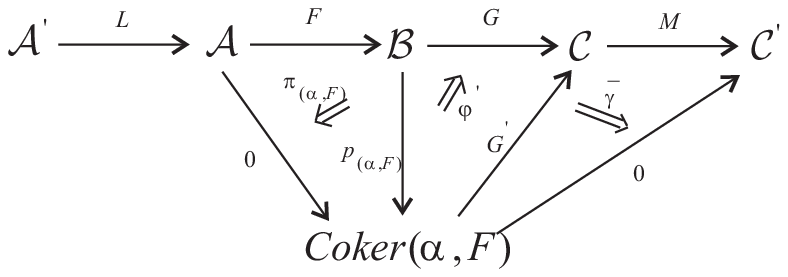}}
\end{center}
We say that the sequence $(L,\alpha,F,\varphi,G,\gamma,M)$ is
relative 2-exact in $\cB$ if the functor $G^{'}$ is faithful and
$\overline{\gamma}$-full.
\end{Def}
\begin{Rek}
The equivalent definition of relative 2-exact can also be given
using relative kernel similar to the symmetric 2-group case in
\cite{11}.
\end{Rek}

In the following, we will omit the composition symbol $\circ$ in our
diagrams.

A complex of $\cR$-2-modules is a diagram in ($\cR$-2-Mod) of the
form
$$
\cA_{\cdot}=\cA_{0}\xrightarrow[]{L_{0}}\cA_{1}\xrightarrow[]
{L_{1}}\cA_{2}\xrightarrow[]{L_{2}}\cdot\cdot\cdot\xrightarrow[]{L_{n-1}}\cA_{n}\xrightarrow[]{L_{n}}\cA_{n+1}
\xrightarrow[]{L_{n+1}}\cA_{n+2}\xrightarrow[]{L_{n+2}}\cdot\cdot\cdot
$$
together with a family of 2-morphisms $\{\alpha_{n}:L_{n+1}\circ
L_{n}\Rightarrow 0\}_{n\geq0}$ such that, for all $n$, the following
diagram commutes
\begin{center}
\scalebox{0.9}[0.85]{\includegraphics{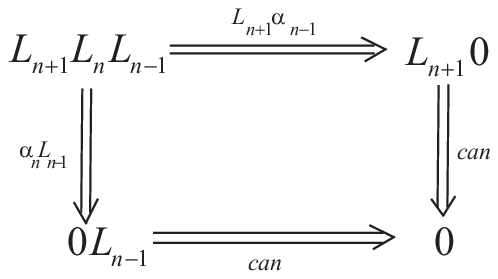}}
\end{center}
We call this complex a 2-cochain complex in ($\cR$-2-Mod).

Consider part of the 2-cochain complex
\begin{center}
\scalebox{0.9}[0.85]{\includegraphics{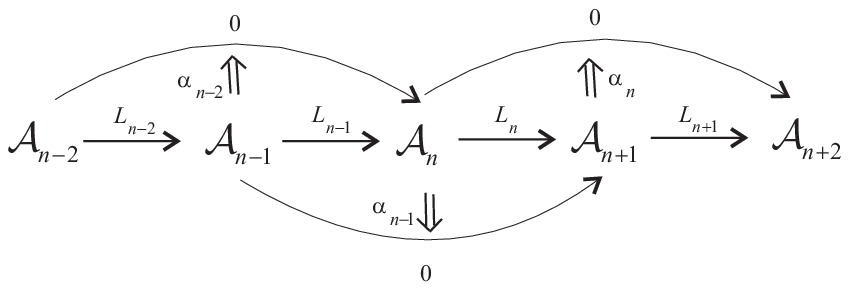}}
\end{center}
Based on the universal properties of relative kernel
$Ker(L_n,\alpha_n)$, we have the following diagram
\begin{center}
\scalebox{0.9}[0.85]{\includegraphics{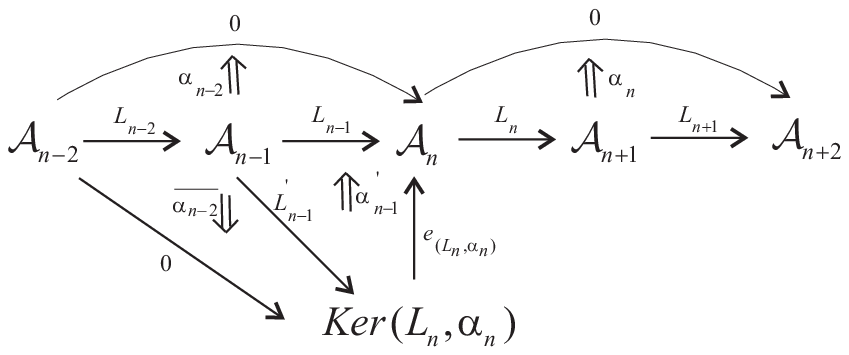}}
\end{center}
The $n$th cohomology $\cR$-2-module $\cH^{n}(\cA_{\cdot})$ of the
2-cochain complex $\cA_{\cdot}$ is defined as the relative cokernel
$Coker(\overline{\alpha_{n-2}},L_{n-1}^{'})$.

Note that, to get $\cH^{0}(\cA_{\cdot})$ and $\cH^{1}(\cA_{\cdot})$,
we have to complete the complex $\cA_{\cdot}$ on the left with the
two zero morphisms and two canonical 2-morphisms\\
$ 0\xrightarrow[]{0}0\xrightarrow
{0}\cA_0\xrightarrow[]{L_0}\cA_{1}\xrightarrow[]{L_1}\cdot\cdot\cdot,
$ can : $L_0\circ0 \Rightarrow 0$ can : $0\circ 0\Rightarrow 0.$

For the convenient application of $\cH^{n}(\cA_{\cdot})$, we will
give an explicit description of it using the construction of
relative kernel and cokernel (similar to \cite{11}).

$\cdot$ An object of $\cH^{n}(\cA_{\cdot})$ is an object of the
relative kernel $Ker(L_{n},\alpha_n)$, that is a pair
$$
(A_n\in obj(\cA_n), a_n:L_n(A_n)\rightarrow 0)
$$
such that $L_{n+1}(a_n)=(\alpha_{n})_{A_n}$;

$\cdot$ A morphism $(A_n,a_n)\rightarrow (A_n^{'},a_n^{'})$ is an
equivalent pair
$$
(X_{n-1}\in obj(\cA_{n-1}),x_{n-1}:A_n\rightarrow
L_{n-1}(X_{n-1})+A_n^{'})
$$
such that the following diagram commutes
\begin{center}
\scalebox{0.9}[0.85]{\includegraphics{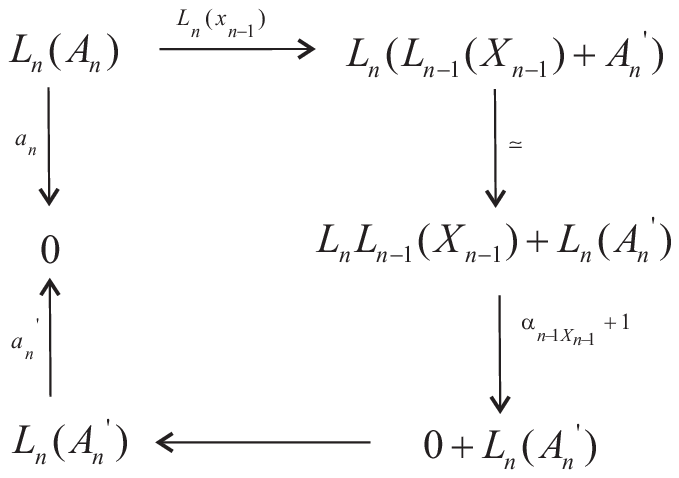}}
\end{center}
Two morphisms
$(X_{n-1},x_{n-1}),(X_{n-1}^{'},x_{n-1}^{'}):(A_n,a_n)\rightarrow
(A_n^{'},a_n^{'})$ are equivalent if there is a pair
$$
(X_{n-2}\in obj(\cA_{n-2}),x_{n-2}:X_{n-1}\rightarrow
L_{n-2}(X_{n-2})+X_{n-1}^{'})
$$
such that the following diagram commutes
\begin{center}
\scalebox{0.9}[0.85]{\includegraphics{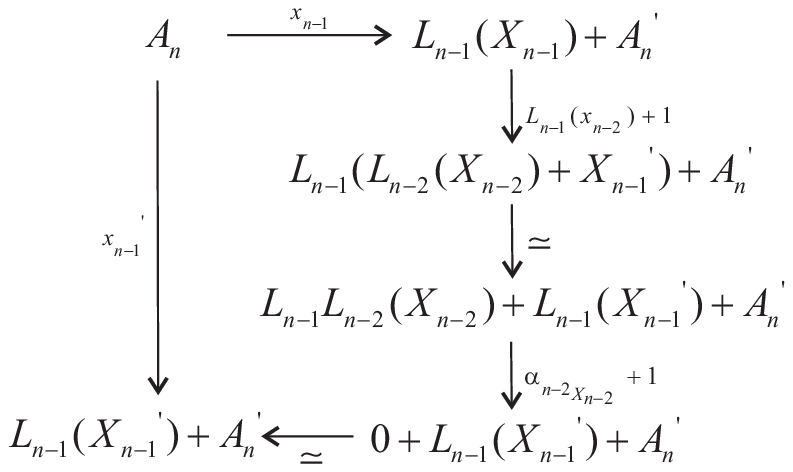}}
\end{center}

A morphism
$(F_{\cdot},\lambda_{\cdot}):\cA_{\cdot}\rightarrow\cB_{\cdot}$ of
two 2-cochain complexes in ($\cR$-2-Mod) is a picture in the
following diagram
\begin{center}
\scalebox{0.9}[0.85]{\includegraphics{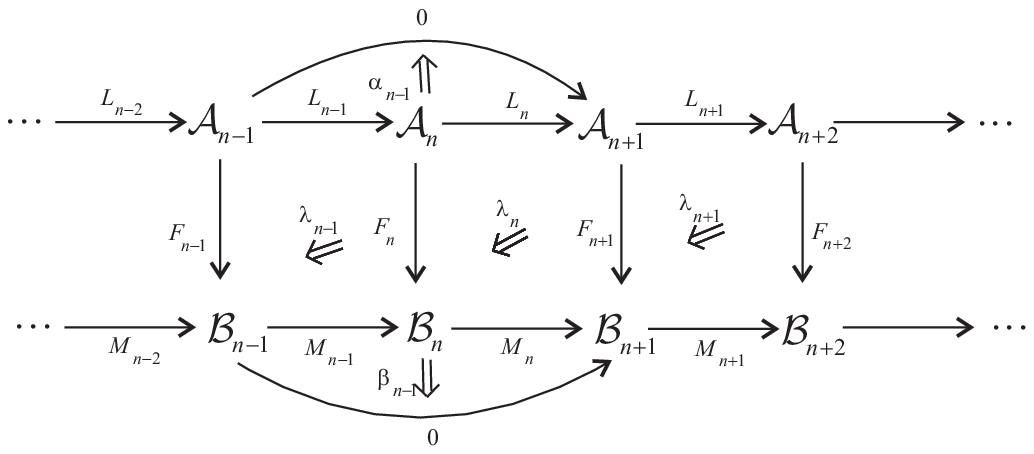}}
\end{center}
where $F_{n}:\cA_n\rightarrow \cB_n$ is 1-morphism in ($\cR$-2-Mod),
$\lambda_{n}:F_{n+1}\circ L_{n}\Rightarrow M_{n}\circ F_n$ is
2-morphism in ($\cR$-2-Mod), for each $n$, making the following
diagram commutative
\begin{center}
\scalebox{0.9}[0.85]{\includegraphics{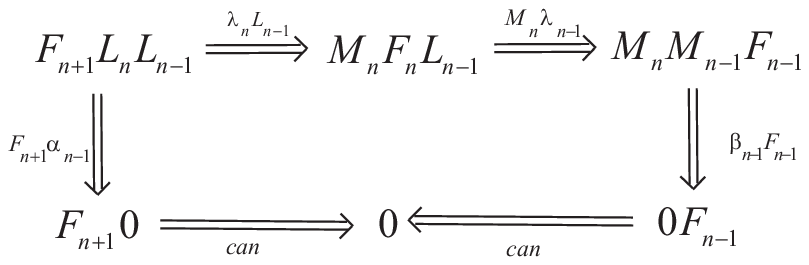}}
\end{center}
Such a morphism induces, for each $n$, a morphism of cohomology
$\cR$-2-modules $\cH^{n}(F_{\cdot}):\cH^{n}(\cA_{\cdot})\rightarrow
\cH^{n}(\cB_{\cdot})$ from the universal properties of relative
kernel and cokernel. $\cH^{n}(F_{\cdot})$ is a morphism of
cohomology symmetric 2-groups, which is given in \cite{11}, together
with a natural isomorphism
$\cH^{n}(F_{\cdot})_{2}:\cH^{n}(F_{\cdot})(r\cdot(A_n,a_n))\rightarrow
r\cdot \cH^{n}(F_{\cdot})(A_n,a_n)$ induced by the
$\cR$-homomorphism $F_n$(\cite{4}).

\begin{Rek}
1. For a complex of $\cR$-2-modules which is relative 2-exact in
each point, the (co)homology $\cR$-2-modules are always zero
$\cR$-2-module(only one object and one morphism)(\cite{2,4}).

2. For morphisms
$\cA_{\cdot}\xrightarrow[]{(F_{\cdot},\lambda_{\cdot})}\cB_{\cdot}\xrightarrow[]{(G_{\cdot},\mu_{\cdot})}\cC_{\cdot}$
of complexes of $\cR$-2-modules, their composite is given by
$(G_{n}\circ F_{n},(\mu_n\circ F_{n})\star (G_{n+1}\circ
\lambda_n))$, for $n\in \mathds{Z}$, where $\star$ is the vertical
composition of 2-morphisms in 2-category(\cite{19,2,4}). Moreover,
$\cH^{n}(G_{\cdot}\circ F_{\cdot})\backsimeq \cH^{n}(G_{\cdot})\circ
\cH^n(F_{\cdot})$ of cohomology $\cR$-2-modules.
\end{Rek}

\begin{Def}
Let $(F_{\cdot},\lambda_{\cdot}),
(G_{\cdot},\mu_{\cdot}):(\cA_{\cdot},L_{\cdot},\alpha_{\cdot})\rightarrow
(\cB_{\cdot},M_{\cdot},\beta_{\cdot})$ be two morphisms of 2-cochain
complexes of $\cR$-2-modules. If there is a family of 1-morphisms
$\{H_{n-1}:\cA_{n}\rightarrow\cB_{n-1}\}_{n\in \mathds{Z}}$ and a
family of 2-morphisms $\{\tau_{n}:F_{n}\Rightarrow M_{n-1}\circ
H_{n-1}+H_{n}\circ L_{n}+G_{n}:\cA_{n}\rightarrow \cB_n\}_{n\in
\mathds{Z}}$ in ($\cR$-2-Mod) satisfying the obvious compatible
conditions, i.e. the following diagram commutes
\begin{center}
\scalebox{0.9}[0.85]{\includegraphics{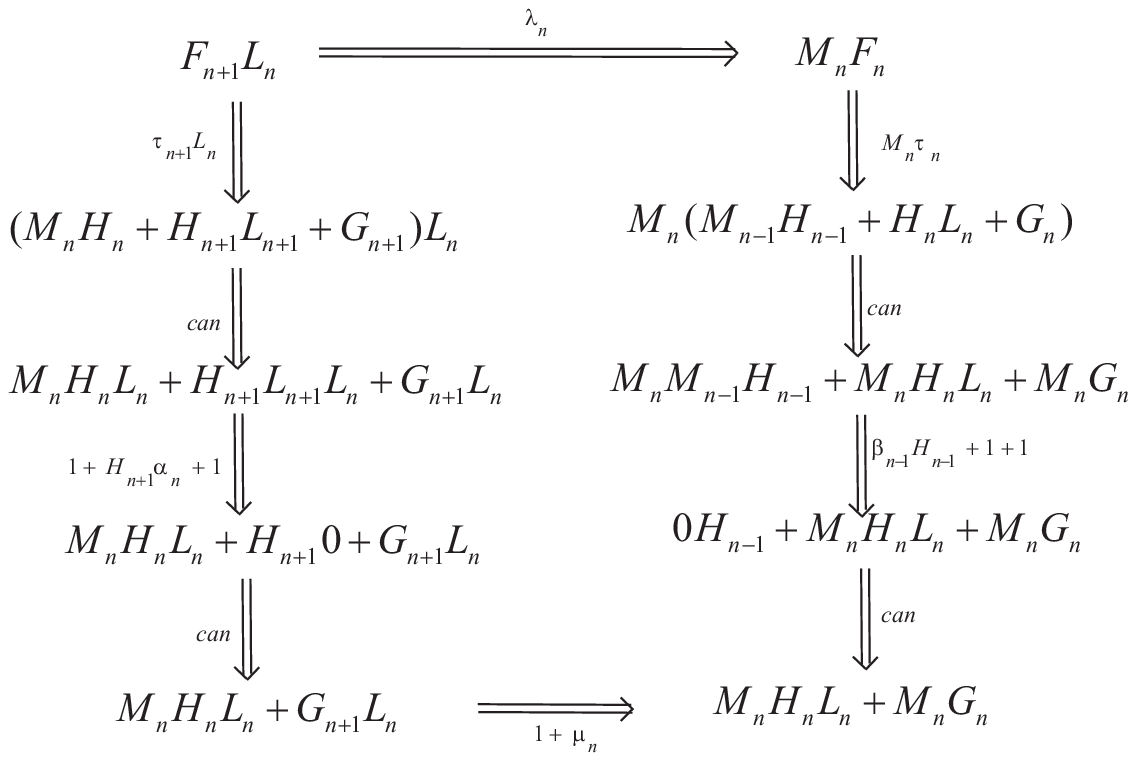}}
\end{center}
We call the above morphisms $(F_{\cdot},\lambda_{\cdot}),
(G_{\cdot},\mu_{\cdot})$ are 2-cochain homotopy.
\end{Def}

\begin{Prop}
Let $(F_{\cdot},\lambda_{\cdot}),
(G_{\cdot},\mu_{\cdot}):(\cA_{\cdot},L_{\cdot},\alpha_{\cdot})\rightarrow
(\cB_{\cdot},M_{\cdot},\beta_{\cdot})$ be two morphisms of 2-cochain
complexes of $\cR$-2-modules. If they are 2-cochain homotopy, there
is an equivalence between induced morphisms $\cH^{n}(F_{\cdot})$ and
$\cH^{n}(G_{\cdot})$.
\end{Prop}
\begin{proof}
In order to prove the equivalence between two morphisms, it will
suffice to construct a 2-morphism
$\varphi_n:\cH^{n}(F_{\cdot})\Rightarrow \cH^{n}(G_{\cdot})$, for
each $n$.

There are induced morphisms
\begin{align*}
&\cH^{n}(F_{\cdot}):\cH^{n}(\cA_{\cdot})\rightarrow \cH^{n}(\cB_{\cdot})\\
&\hspace{1.5cm}(A_n,a_n)\mapsto(F_n (A_n),b_n),\\
&\hspace{0.8cm}[X_{n-1},x_{n-1}]\mapsto
[F_{n-1}(X_{n-1}),\overline{x_{n-1}}]
\end{align*}
and
\begin{align*}
&\cH^{n}(G_{\cdot}):\cH^{n}(\cA_{\cdot})\rightarrow \cH^{n}(\cB_{\cdot})\\
&\hspace{1.5cm}(A_n,a_n)\mapsto(G_n (A_n),\overline{b_n}),\\
&\hspace{0.8cm}[X_{n-1},x_{n-1}]\mapsto
[G_{n-1}(X_{n-1}),\overline{x_{n-1}}^{'}]
\end{align*}
where $b_n$ is the composition
$M_{n}F_{n}(A_n)\xrightarrow[]{(\lambda_{n})_{A_n}^{-1}}F_{n+1}L_{n}(A_n)\xrightarrow[]{F_{n+1}(a_n)}
F_{n+1}(0)\backsimeq 0$, $\overline{x_{n-1}}$ is the composition
$F_{n}(A_n)\xrightarrow[]{F_{n}(x_{n-1})}F_{n}(L_{n-1}(X_{n-1})+A_{n}^{'})\backsimeq
F_{n}L_{n-1}(X_{n-1})+F_{n}(A_{n}^{'})\xrightarrow[]{(\lambda_{n-1})_{X_{n-1}}+1}M_{n-1}F_{n-1}(X_{n-1})+F_{n}(A_{n}^{'}),$
$\overline{b_n}$ and $\overline{x_{n-1}}^{'}$ can be  given in
similar ways. For any object $(A_n,a_n)$ of $\cH^{n}(\cA_{\cdot})$,
let $Y_{n-1}=H_{n-1}(A_n)$. Consider the following composition
morphism $y_{n-1}:F_n(A_n)\xrightarrow[]{(\tau_{n})_{A_n}}
(M_{n-1}\circ H_{n-1}+H_{n}\circ
L_n+G_n)(A_n)\xrightarrow[]{1+H_{n}(a_n)+1}M_{n-1}(H_{n-1}(A_n))+H_{n}(0)+G_n(A_n)\backsimeq
M_{n-1}(Y_{n-1})+0+G_n(A_n)\backsimeq M_{n-1}(Y_{n-1})+G_n(A_n)$ of
$\cB_{n}$. We get a morphism $[Y_{n-1}\in obj(\cB_{n-1}),
y_{n-1}:F_{n}(A_n)\rightarrow
M_{n-1}(Y_{n-1})+G_{n}(A_n)]:(F_{n}(A_n),b_n)\rightarrow
(G_{n}(A_n),\overline{b_{n}})$ of $\cH^{n}(\cB_{\cdot})$ such that
the following diagram commutes
\begin{center}
\scalebox{0.9}[0.85]{\includegraphics{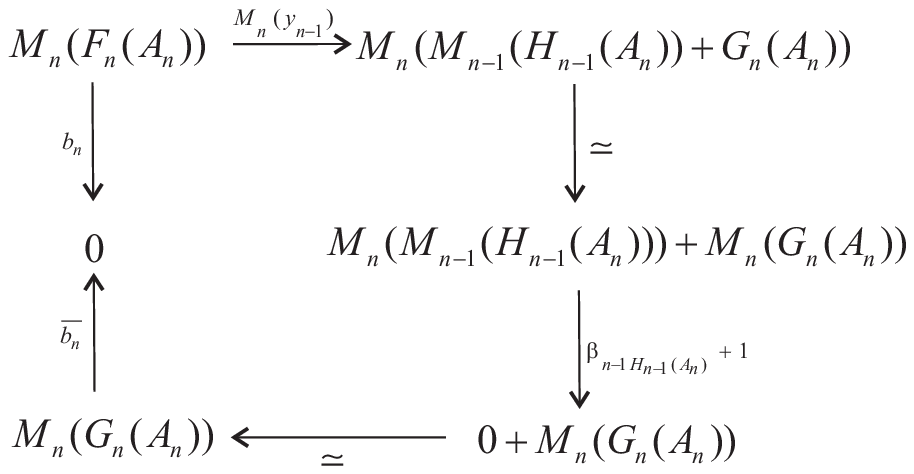}}
\end{center}
From the compatible condition of $(\tau_{n})_{n\in\mathds{Z}}$, we
have the following commutative diagram
\begin{center}
\scalebox{0.7}[0.7]{\includegraphics{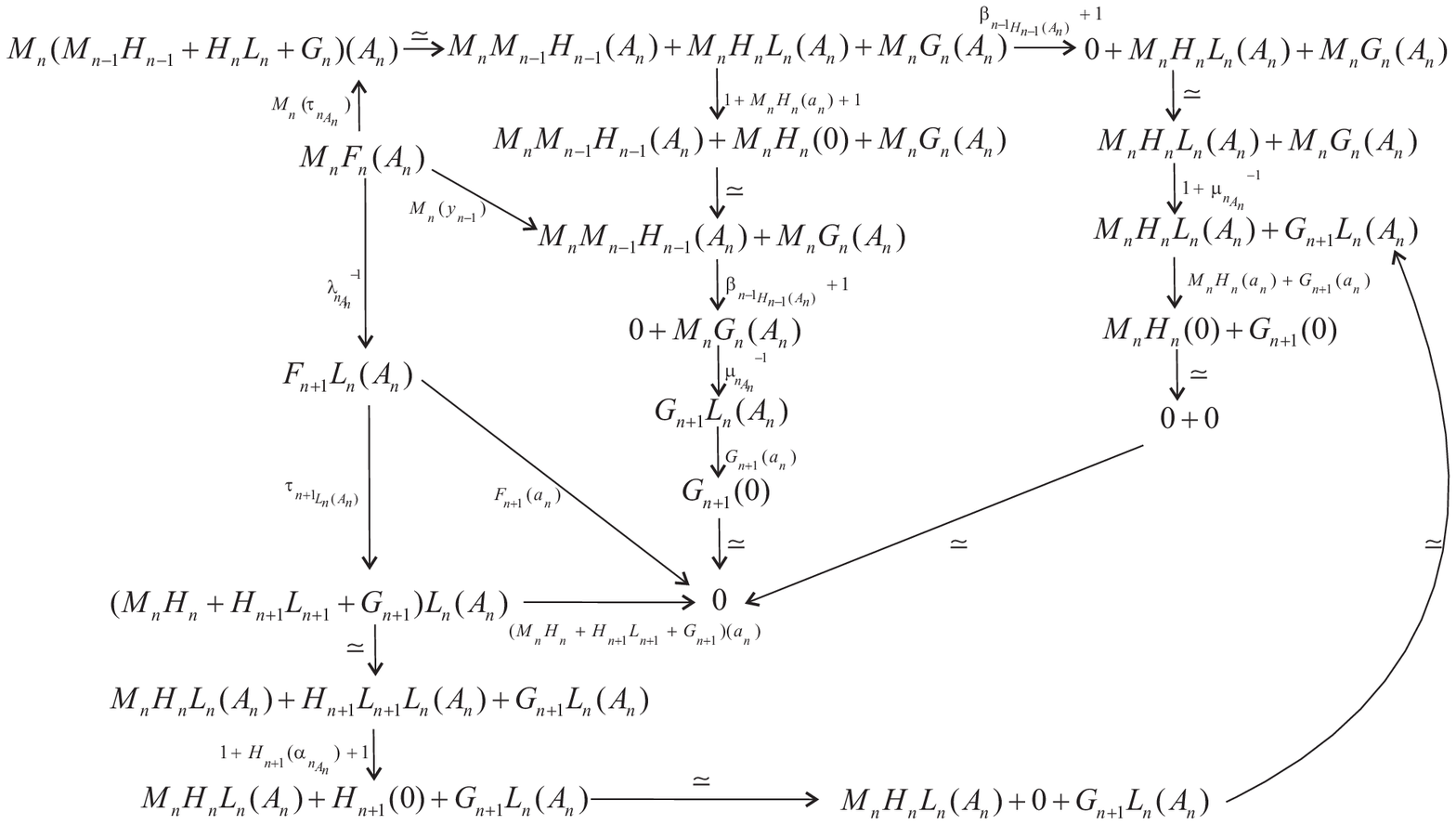}}
\end{center}
So $[Y_{n-1},y_{n-1}]$ is a morphism in $\cH^{n}(\cB_{\cdot})$, then
we can define a 2-morphism $\varphi_{n}$. For any morphism
$[X_{n-1},x_{n-1}]:(A_n,a_n)\rightarrow (A_{n}^{'},a_{n}^{'})$ in
$\cH^{n}(\cA_{\cdot})$, where $X_{n-1}\in obj(\cA_{n+1}),\
x_{n-1}:A_n\rightarrow L_{n-1}(X_{n-1})+A_{n}^{'}$ satisfying the
following commutative diagram
\begin{center}
\scalebox{0.9}[0.85]{\includegraphics{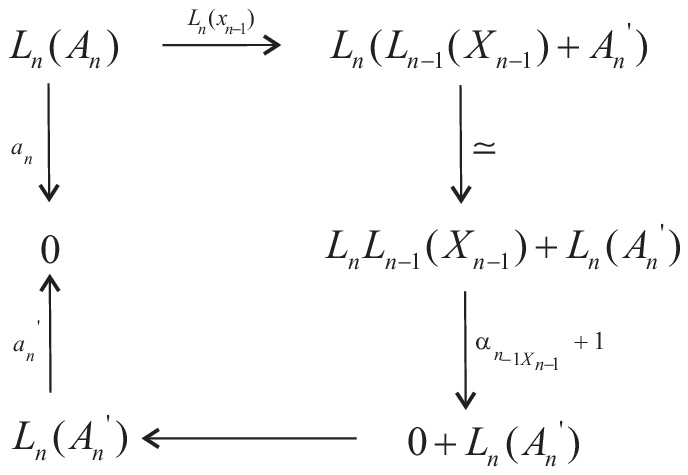}}
\end{center}
$\cH^{n}(F_{\cdot})[X_{n-1},x_{n-1}]=[F_{n-1}(X_{n-1}),\overline{x_{n-1}}]$,
$\cH^{n}(G_{\cdot})[X_{n-1},x_{n-1}]=[G_{n-1}(X_{n-1}),\overline{x_{n-1}}^{'}]$,
where $\overline{x_{n-1}}$ and $\overline{x_{n-1}}^{'}$ are the
following composition morphisms
$\overline{x_{n-1}}:F_{n}(A_n)\xrightarrow[]{F_{n}(x_{n-1})}F_{n}(L_{n-1}(X_{n-1})+A_{n}^{'})\backsimeq
F_{n}L_{n-1}(X_{n-1})+F_{n}(A_{n}^{'})\xrightarrow[]{(\lambda_{n-1})_{X_{n-1}}+1}M_{n-1}F_{n-1}(X_{n-1})+F_{n}(A_{n}^{'})$,
$\overline{x_{n+1}}^{'}:G_{n}(A_n)\xrightarrow[]{G_{n}(x_{n-1})}G_{n}(L_{n-1}(X_{n-1})+A_{n}^{'})\backsimeq
G_{n}L_{n-1}(X_{n-1})+G_{n}(A_{n}^{'})\xrightarrow[]{(\mu_{n-1})_{X_{n-1}}+1}M_{n-1}G_{n-1}(X_{n-1})+G_{n}(A_{n}^{'}).$

Then we have the following commutative diagram
\begin{center}
\scalebox{0.9}[0.85]{\includegraphics{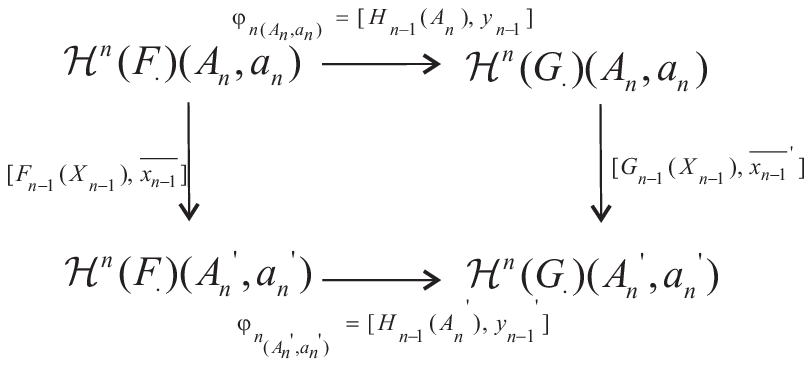}}
\end{center}
There exist $[Y_{n-2}\triangleq
H_{n-2}(X_{n-1}),y_{n-2}]:([H_{n-1}(A_{n}^{'}),y_{n-1}^{'}]\circ
[F_{n-1}(X_{n-1}),\overline{x_{n-1}}])\rightarrow
[G_{n-1}(X_{n-1}),\overline{x_{n-1}}^{'}]\circ[H_{n-1}(A_{n}),y_{n-1}]$
induced by $\tau_{\cdot}.$ In fact,
$((H_{n-1}(A_{n}^{'}),y_{n-1}^{'})\circ
F_{n-1}(X_{n-1}),\overline{x_{n-1}}))=[F_{n-1}(X_{n-1})+H_{n-1}(A_{n}^{'}),(1+y_{n-1}^{'})\circ
\overline{x_{n-1}}^{'}]$,
$[G_{n-1}(X_{n-1}),\overline{x_{n-1}}]\circ[H_{n-1}(A_{n}),y_{n-1}]=[H_{n-1}(A_n)+G_{n-1}(X_{n-1}),(1+\overline{x_{n-1}}^{'})\circ
y_{n-1}]$ from the composition of morphisms in relative cokernel, so
$y_{n-2}$ is the composition morphism
$F_{n-1}(X_{n-1})+H_{n-1}(A_{n}^{'})\xrightarrow[]{(\tau_{n-1})_{X_{n-1}}+1}(M_{n-2}H_{n-2}+H_{n-1}L_{n-1}+G_{n-1})(X_{n-1})+H_{n-1}(A_{n}^{'})\backsimeq
M_{n-2}H_{n-2}(X_{n-1})+H_{n-1}L_{n-1}(X_{n-1})+G_{n-1}(X_{n-1})+H_{n-1}(A_{n}^{'})\backsimeq
M_{n-2}H_{n-2}(X_{n-1})+G_{n-1}(X_{n-1})+H_{n-1}(L_{n-1}(X_{n-1}+A_{n}^{'})\xrightarrow[]{1+1+H_{n-1}(x_{n-1}^{-1})}
M_{n-2}H_{n-2}(X_{n-1})+G_{n-1}(X_{n-1})+H_{n-1}(A_{n})\backsimeq
M_{n-2}H_{n-2}(X_{n-1})+H_{n-1}(A_{n})+G_{n-1}(X_{n-1})$. Moreover
the morphism $[Y_{n-2},y_{n-2}]$ makes the following diagram commute
\begin{center}
\scalebox{0.8}[0.8]{\includegraphics{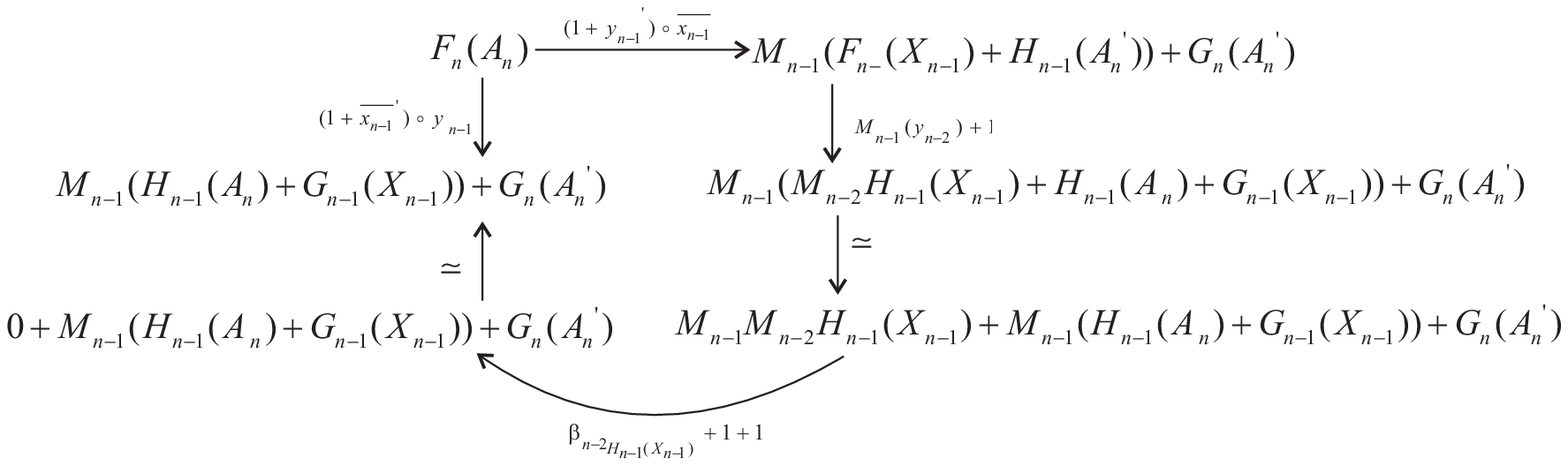}}
\end{center}
from the following several commutative diagrams
\begin{center}
\scalebox{0.8}[0.8]{\includegraphics{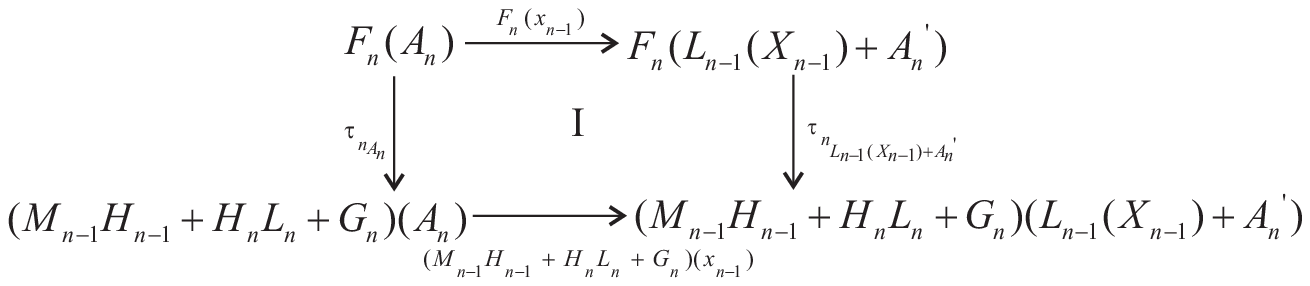}}
\end{center}
\begin{center}
\scalebox{0.8}[0.8]{\includegraphics{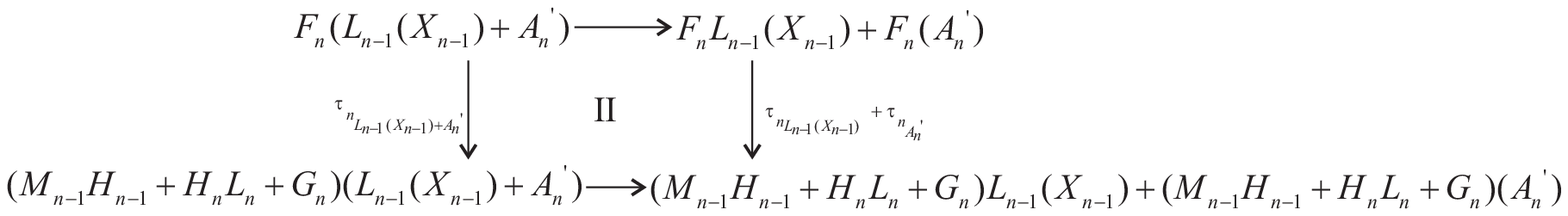}}
\end{center}
\begin{center}
\scalebox{0.8}[0.8]{\includegraphics{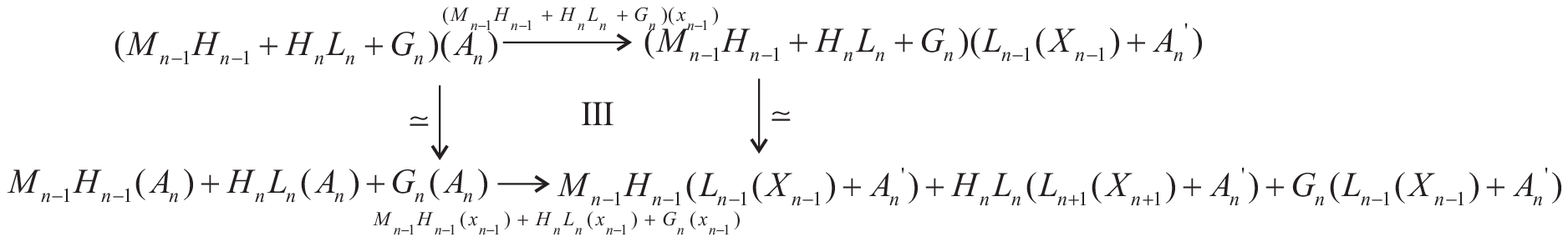}}
\end{center}
\begin{center}
\scalebox{0.7}[0.7]{\includegraphics{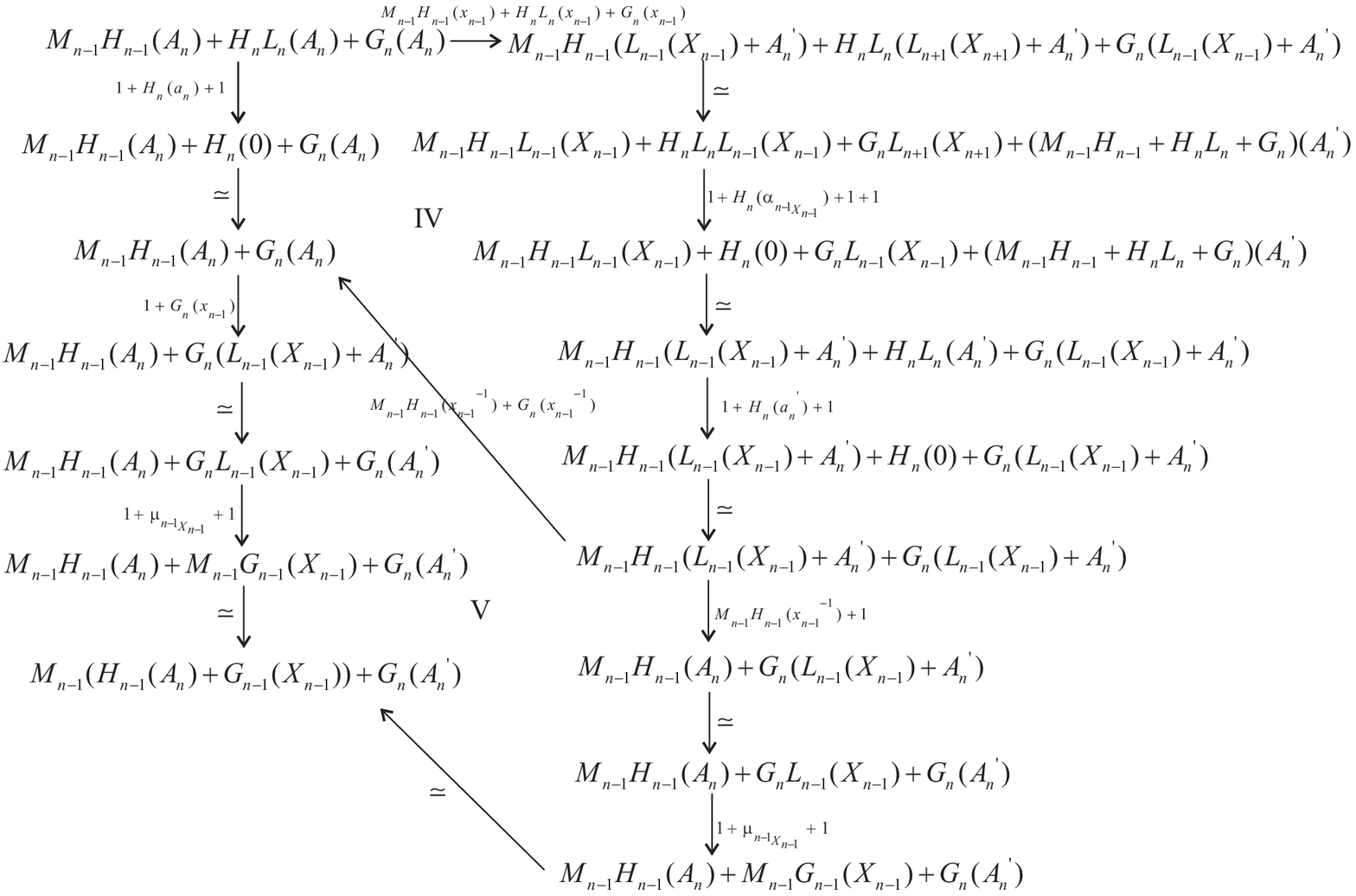}}
\end{center}
\begin{center}
\scalebox{0.7}[0.7]{\includegraphics{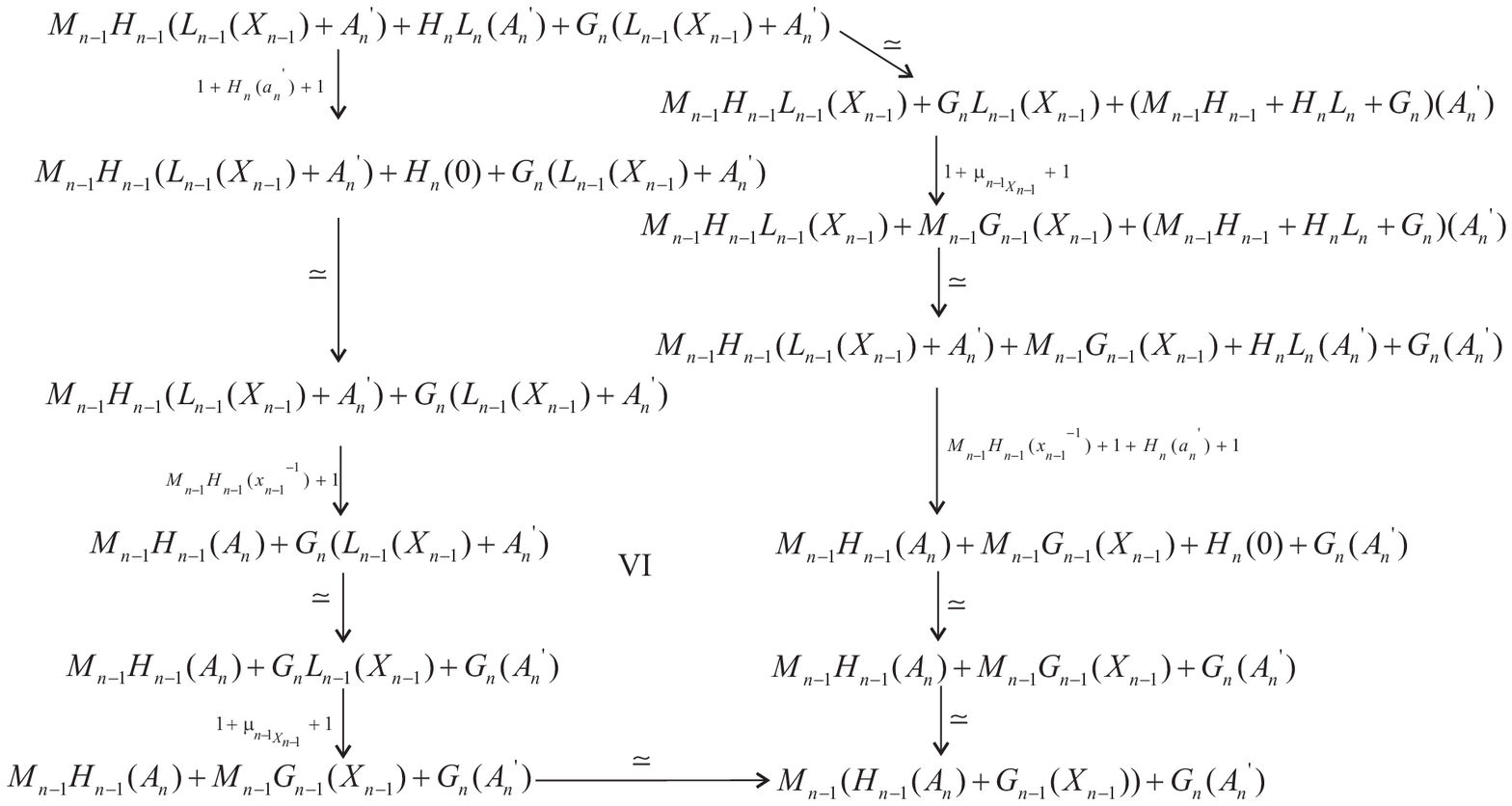}}
\end{center}
\begin{center}
\scalebox{0.7}[0.7]{\includegraphics{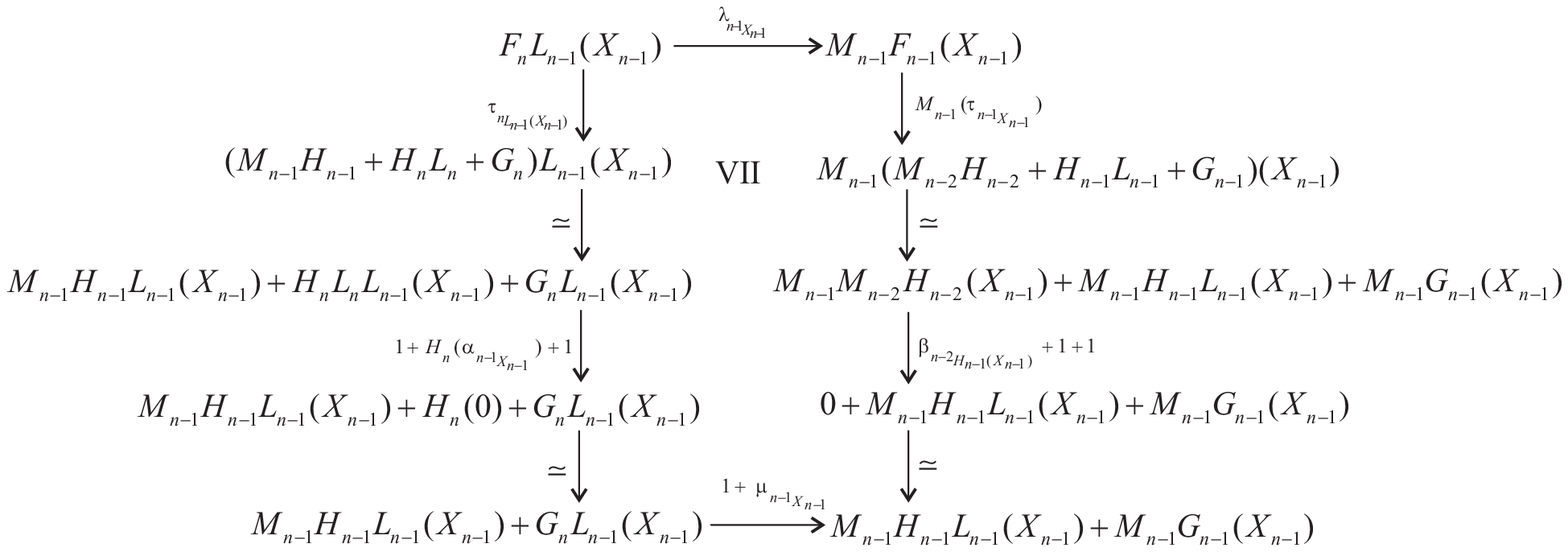}}
\end{center}
\begin{center}
\scalebox{0.7}[0.7]{\includegraphics{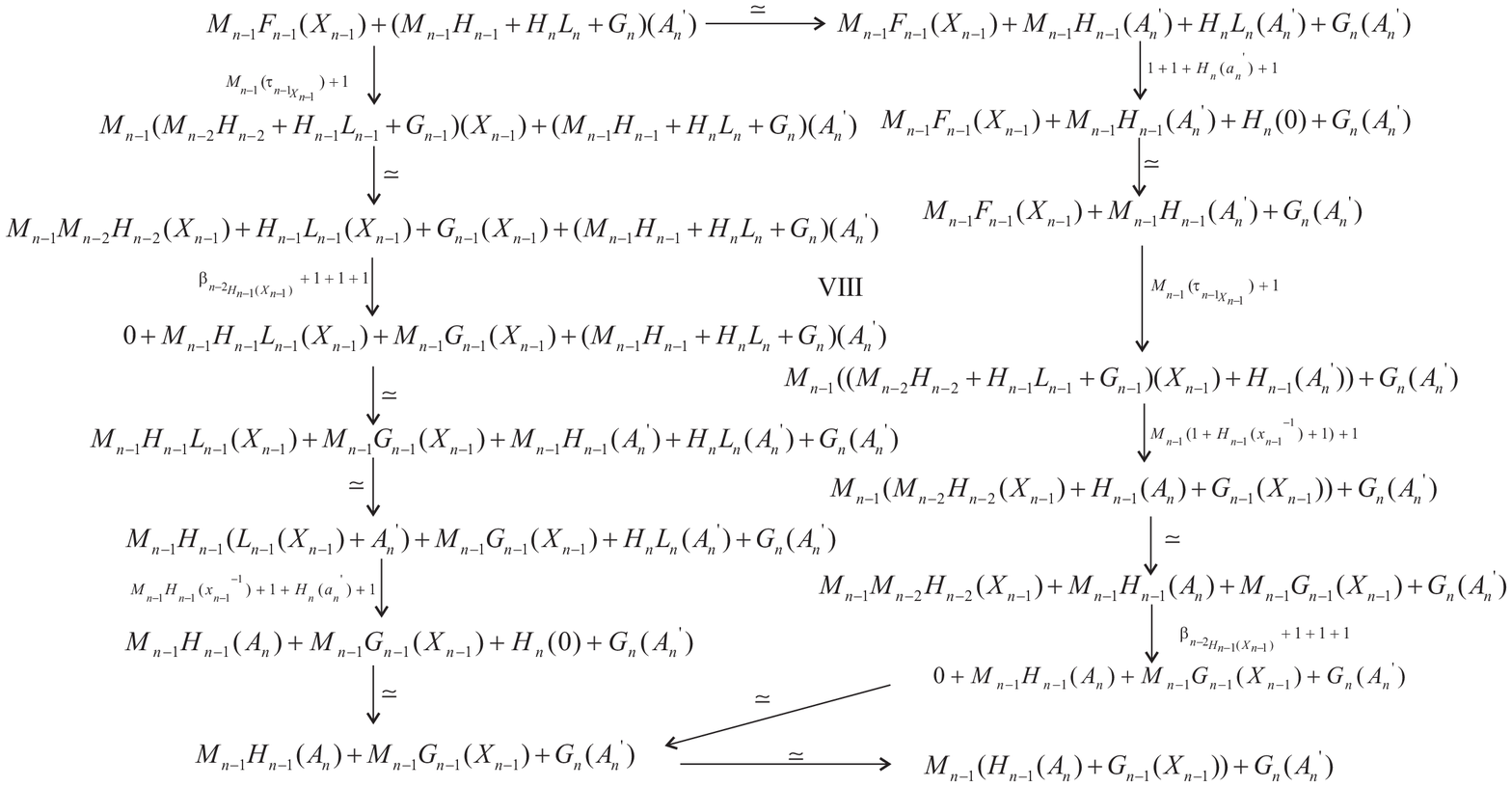}}
\end{center}
where I is commutative because $\tau_n$ is a natural transformation.
II, III, VIII follow from the properties of symmetric monoidal
functors. IV follows from operation of $H_{n}$ on the commutative
diagram of $[X_{n-1},x_{n-1}]$. V and VI follow from the properties
of symmetric 2-groups. VII follows from the commtutative diagram of
$\tau_{n}$.

For any two objects $(A_n,a_n),\ (A_n^{'},a_n^{'})$ of
$\cH^{n}(\cA_{\cdot})$,
$(A_n,a_n)+(A_n^{'},a_n^{'})=(A_{n}+A_{n}^{'},a_{n}+a_n^{'})$ with
$L_{n+1}(a_n+a_{n}^{'})=(\alpha_{n})_{A_n+A_n^{'}}$, we have the
following commutative diagram
\begin{center}
\scalebox{0.9}[0.85]{\includegraphics{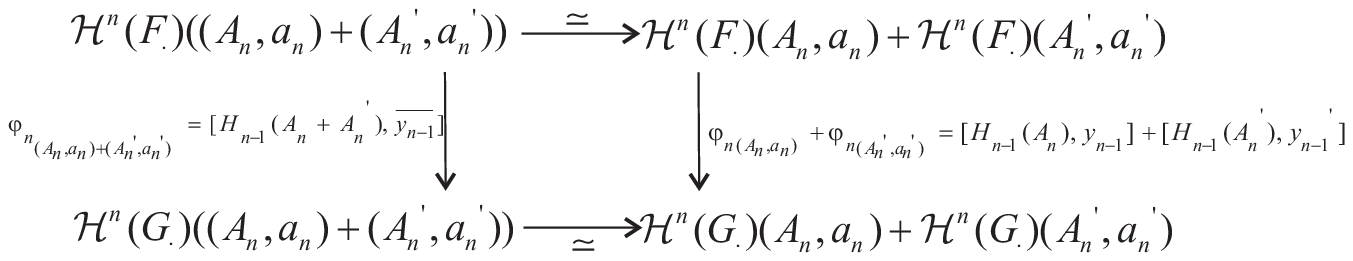}}
\end{center}
where $\overline{y_{n-1}}, y_{n-1},y_{n-1}^{'}$ are induced by
$\tau_{n}$ as above. In fact,
$[H_{n-1}(A_n),y_{n-1}]+[H_{n-1}(A_{n}^{'}),y_{n-1}^{'}]=[H_{n-1}(A_n)+H_n(A_{n}^{'}),y_{n-1}+y_{n-1}^{'}]$,
so
$(\varphi_{n})_{(A_n,a_n)+(A_n^{'},a_{n}^{'})}=(\varphi_{n})_{(A_n,a_n)}+(\varphi_{n})_{(A_n^{'},a_n^{'})}$,
then the above diagram commutes.

Moreover, for any $r\in obj(\cR),\ (A_n,a_n)\in
obj(\cH^{n}(\cA_{\cdot}))$, there is a commutative diagram
\begin{center}
\scalebox{0.9}[0.85]{\includegraphics{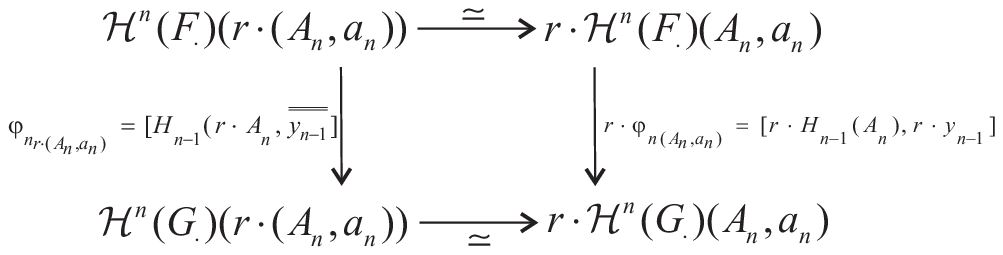}}
\end{center}
following from $H_{n-1}$ is an $\cR$-homomorphism.

Then $\varphi_n$ is a 2-morphism in ($\cR$-2-Mod), for each $n$.
\end{proof}

From the properties of 2-functors and natural proof, we get
\begin{Lem}Let $\cR,\cS$ be 2-rings, $(F_{\cdot},\lambda_{\cdot}),
(G_{\cdot},\mu_{\cdot}):(\cA_{\cdot},L_{\cdot},\alpha_{\cdot})\rightarrow
(\cB_{\cdot},M_{\cdot},\beta_{\cdot})$ be two 2-cochain homotopic
morphisms of 2-cochain complexes of $\cR$-2-modules and
$T:(\cR$-2-Mod)$\rightarrow(\cS$-2-Mod) be a 2-functor. Then
$T(F_{\cdot})$ is 2-cochain homotopic to $T(G_{\cdot})$  in
($\cS$-2-Mod).
\end{Lem}

\section{Injective Resolution of $\cR$-2-Module}

In this section we will give the construction of injective
resolution of any $\cR$-2-module.

\begin{Def} Let $\cR$ be a 2-ring, and $\cA$ be an $\cR$-2-module.
An injective resolution of $\cA$ in ($\cR$-2-Mod) is a 2-cochain
complex of $\cR$-2-modules which is relative 2-exact in each point
as in the following diagram
\begin{center}
\scalebox{0.9}[0.85]{\includegraphics{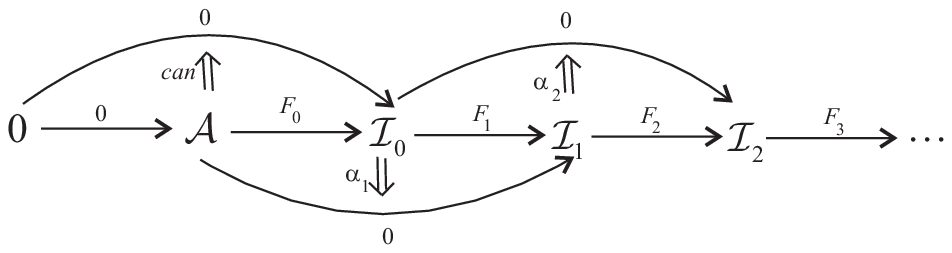}}
\end{center}
with  injective objects $\cI_{n}(n\geq 0)$in ($\cR$-2-Mod). i.e. the
above complex is relative 2-exact in $\cA$ and each $\cI_{n}$, for
$n\geq 0$.
\end{Def}

\begin{Prop}
Every $\cR$-2-module $\cA$ has an injective resolution in
($\cR$-2-Mod).
\end{Prop}
\begin{proof}
We will construct the injective resolution of $\cA$ using the
relative cokernel.

For $\cA$, there is a faithful morphism $F_0:\cA\rightarrow\cI_{0}$,
with $\cI_{0}$ injective in ($\cR$-2-Mod)(\cite{20}). Then we get a
sequence as follows
\begin{center}
\scalebox{0.9}[0.85]{\includegraphics{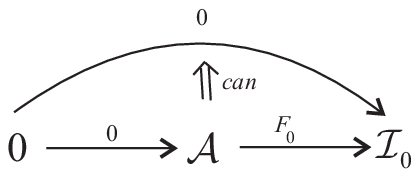}} {\footnotesize S.1.}
\end{center}
where $0:0\rightarrow\cA$ is the zero morphism\cite{2,4,6}in
($\cR$-2-Mod), $0$ is the $\cR$-2-module with only one object and
one morphism., $can$ is the canonical 2-morphism in ($\cR$-2-Mod),
which is given by the identity morphism of only one object of $0$.

From the existence of the relative cokernel in ($\cR$-2-Mod), we
have the relative cokernel $(Coker(can,F_0),\ p_{(can,F_0)},\
\pi_{(can,F_0)})$
 of the sequence S.1, which is in fact the general cokernel
 $(CokerF_0,p_{F_0},\pi_{F_0})$ \cite{2,28}. For the $\cR$-2-module $CokerF_0$, there exists a faithful morphism
 $\overline{F_1}:CokerF_0\rightarrow \cI_{1}$, with $\cI_{1}$ injective in ($\cR$-2-Mod)(\cite{20}).
 Let $F_{1}=\overline{F_{1}}\circ p_{F_{0}}:\cI_{0}\rightarrow\cI_{1}$. Then we get the following
 sequence
\begin{center}
\scalebox{0.9}[0.85]{\includegraphics{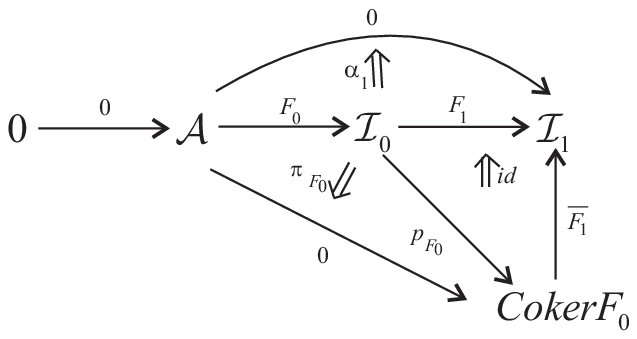}}
\end{center}
where $\alpha_{1}$ is the composition $F_{1}\circ F_{0}=
\overline{F_1}\circ p_{F_0}\circ F_{0}\Rightarrow
\overline{F_1}\circ0\Rightarrow 0$ and is compatible with $can$.

Consider the above sequence, there exists the relative cokernel
$(Coker(\alpha_1,F_1),\\p_{(\alpha_1,F_1)},\pi_{(\alpha_1,F_1})$ in
($\cR$-2-Mod). For the $\cR$-2-module $Coker(\alpha_1,F_1)$, there
is a faithful morphism
$\overline{F_{2}}:Coker(\alpha_1,F_1)\rightarrow \cI_{2}$, with
$\cI_2$ injective in ($\cR$-2-Mod)(\cite{20}). Let
$F_2=\overline{F_2}\circ p_(\alpha_1,F_1):\cI_1\rightarrow\cI_2$.
Then we get a sequence
\begin{center}
\scalebox{0.9}[0.85]{\includegraphics{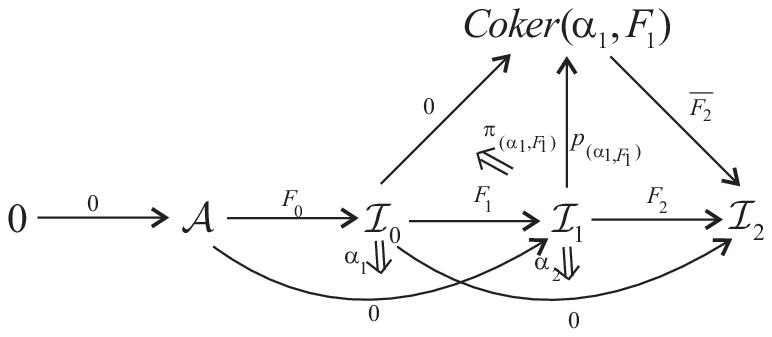}}
\end{center}
where $\alpha_2$ is the composition $F_{2}\circ
F_{1}=\overline{F_{2}}\circ p_{(\alpha_1,F_1)}\circ F_{1}\Rightarrow
\overline{F_2}\circ0\Rightarrow 0$ and is compatible with
$\alpha_1$.

Using the same method, we get a 2-cochain complex of $\cR$-2-modules
\begin{center}
\scalebox{0.9}[0.85]{\includegraphics{p12.eps}}
\end{center}
Next, we will check that this complex is relative 2-exact in each
point.

Since $F_0$ is faithful, this complex is relative 2-exact in $\cA$.

From the cancellation property of $p_{F_0}$, there exists
$\overline{\alpha_2}:F_2\circ \overline{F_1}\Rightarrow 0$ defined
by $\overline{(\alpha_{2})_{y}}\triangleq (\alpha_{2})_{y}:F_2\circ
\overline{F_1}(y)\rightarrow 0,\ \forall y\in obj(\cI_{0})$. And
$\overline{F_{1}}:CokerF_0\rightarrow\cI_{1}$ is in fact
$\overline{F_{1}}(x)=F_{1}(x)$, since $p_{F_0}(x)=x$. For any
$X,Y\in obj(CokerF_0)$ and the morphism
$g:\overline{F_{1}}(X)\rightarrow \overline{F_{1}}(Y)$ of $\cI_1$,
such that the following diagram commutes

\begin{center}
\scalebox{0.9}[0.85]{\includegraphics{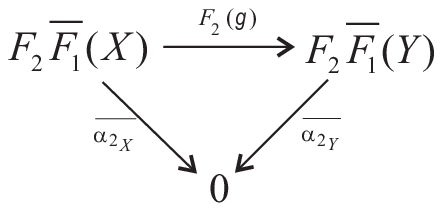}}
\end{center}
From $F_1=\overline{F_1}\circ p_{F_0}$ and the existence of
$\overline{\alpha_2}$, we get a morphism $g:F_{1}(X)\rightarrow
F_{1}(Y)$ in $\cI_1$ and the commutative diagram

\begin{center}
\scalebox{0.9}[0.85]{\includegraphics{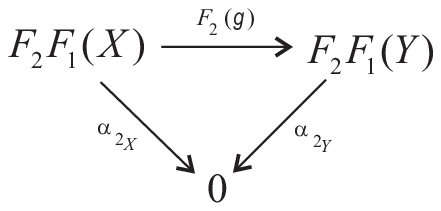}}
\end{center}

Moreover, from the definition of $\alpha_2=\overline{F_2}\circ
\pi_{(\alpha_1,F_1)}$, we have
$(\alpha_2)_X=\overline{F_2}((\pi_{(\alpha_1,F_1)})_{X})=\overline{F_2}([X,1_{F_{1}(X)}])$
and
$(\alpha_2)_Y=\overline{F_2}((\pi_{(\alpha_1,F_1)})_{Y})=\overline{F_2}([Y,1_{F_{1}(Y)}])$,
and then $(\alpha_2)_{Y}^{-1}=\overline{F_1}([Y^{*},1_0])$, where
$Y^{*}$ is the inverse of $Y$ in $\cI_0$. Then, there is
$\overline{F_{2}}p_{(\alpha_1,F_1)}(g)=(\alpha_{2})_{Y}^{-1}\circ
(\alpha_{2})_{X}$, i.e.
$\overline{F_2}([0,g])=\overline{F_2}([Y^{*},1_0]\circ
[X,1_{F_{1}(X)}])=\overline{F_2}([X+Y^{*},1_{F_{1}(X)}])$ from
$\overline{F_2}p_{(\alpha_1,F_1)}F_{1}(X)$ to
$\overline{F_2}p_{(\alpha_1,F_1)}F_{1}(Y)$. Then we have
$[0,g]=[X+Y^{*},1_{F_{1}(X)}]:F_{1}(X)\rightarrow F_{1}(Y)$ from the
faithful morphism $\overline{F_2}$, so there exist $B\in obj(\cA)$
and a morphism $f:0\rightarrow F_{0}(B)+X+Y^{*}$ in $\cI_0$ such
that the following diagram commutes
\begin{center}
\scalebox{0.9}[0.85]{\includegraphics{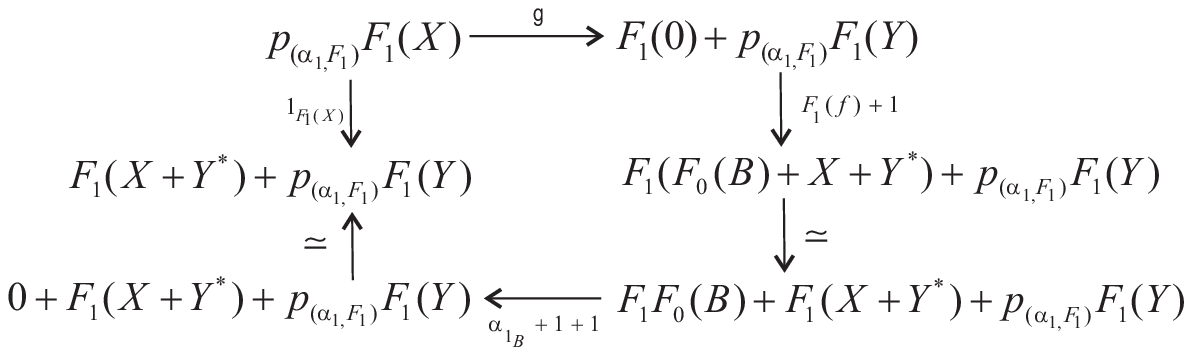}}
\end{center}

Consider the pair $(B,f)$, which is the morphism $[B,f]:0\rightarrow
X+Y^{*}$ in $CokerF_0$, it can be written as a composition
$$
[B,f]=([B,1_{F_{0}(B)}]+1_{X+Y^{*}})\circ [0,f].
$$
Then there is
$\overline{F_{1}}([B,f])=\overline{F_1}(([B,1_{F_{0}(B)}]+1_{X+Y^{*}})\circ
[0,f])=\overline{F_1}(([B,1_{F_{0}(B)}]+1_{X+Y^{*}}))\circ\overline{F_1}(
[0,f])=(\overline{F_1}(([B,1_{F_{0}(B)}])+1_{\overline{F_1}(X+Y^{*}}))\circ\overline{F_1}(
[0,f])=(\overline{F_1}((\pi_{F_0)_{B}})+1_{\overline{F_1}(X+Y^{*}}))\circ\overline{F_1}(
[0,f])=((\alpha_{1})_{B}+1_{\overline{F_1}(X+Y^{*}}))\circ
F_{1}(f).$

Also, from the morphism $[B,f]:0\rightarrow X+Y^{*}$ in $CokerF_0$
and the basic calculations, we have
$[B,f]^{-1}=[B^{*},1_{F_{0}(B)^{-1}}+f^{-1}]:X+Y^{*}\rightarrow 0,$
 $[B,f]+1_{Y}:Y\rightarrow X$ and $[B^*,\overline{f}]:=([B,f]+1_{Y})^{-1}=[B^{*},1_{F_{0}(B)^{-1}}+f^{-1}+1_{Y}]:X\rightarrow
 Y$ in $CokerF_0$.
From the above commutative diagram, we have
$((\alpha_1)_{B}+1_{F_{1}(X+Y^{*})}+1_{F_{1}(Y)})\circ
(F_1(f)+1_{F_{1}(Y)})\circ g=1_{F_{1}(X)}$,  after calculations,
there is $g=\overline{F_1}([B^{*},\overline{f}])$. Then, there exist
a morphism $[B^{*}\in obj(\cA),\overline{f}]:X\rightarrow Y$ in
$CokerF_0$ such that $\overline{F_1}([B^{*},\overline{f}])=g$, i.e.
the faithful morphism $\overline{F_1}$ is
$\overline{\alpha_2}$-full. Then the above complex is relative
2-exact in $\cI_0$.

Using the same method, we can prove the above 2-cochain complex is
relative 2-exact in each point.
\end{proof}

\begin{Thm}
Let $\cR$ be a 2-ring, $(L_{\cdot}:
\cA\rightarrow\cI_{\cdot},\alpha_{\cdot})$ be an injective
resolution of $\cR$-2-module $\cA$, and $F:\cA\rightarrow \cB$ be a
1-morphism in ($\cR$-2-Mod). Then for any injective resolution
$(M_{\cdot}:\cB\rightarrow\cJ_{\cdot},\beta_{\cdot})$ of $\cB$,
there is a morphism $F_{\cdot}:\cI_{\cdot}\rightarrow \cJ_{\cdot}$
of 2-cochain complexes in ($\cR$-2-Mod) together with the family of
2-morphisms $\{\varepsilon_{n}:F_{n}\circ L_{n}\Rightarrow
M_{n}\circ F_{n-1}\}_{n\geq 0}$(where $F_{-1}=F$) as in the
following diagram
\begin{center}
\scalebox{0.9}[0.85]{\includegraphics{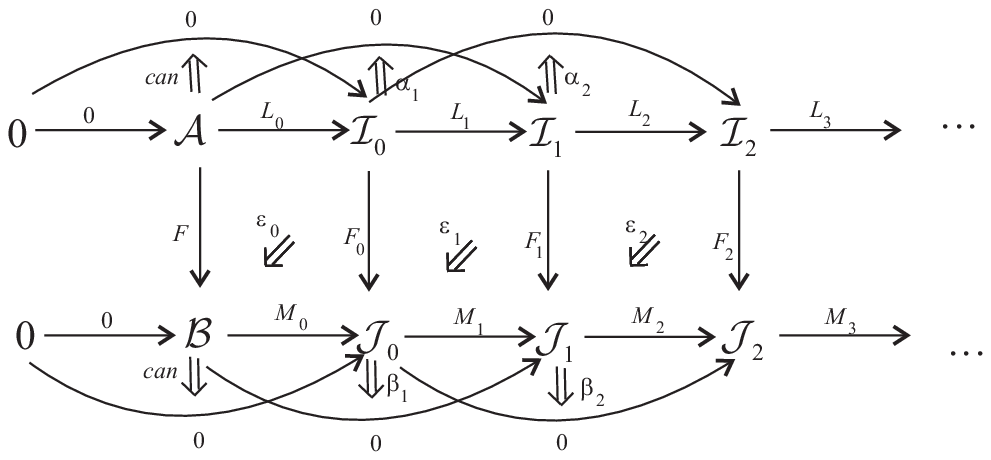}}
\end{center}
If there is another morphism between injective resolutions, they are
2-cochain homotopy.
\end{Thm}
\begin{proof}
The existence of $F_{0}:\cI_{0}\rightarrow \cJ_{0}$: Since $L_{0}$
is faithful and $\cJ_0$ is an injective object in ($\cR$-2-Mod),
there exist 1-morphism $F_0:\cI_{0}\rightarrow \cJ_{0}$ and
2-morphism $\varepsilon_{0}:F_{0}\circ L_{0}\Rightarrow M_0\circ F$
as follows
\begin{center}
\scalebox{0.9}[0.85]{\includegraphics{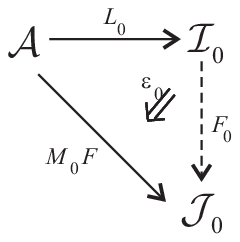}}
\end{center}

Consider the morphism $F_0:\cI_{0}\rightarrow \cJ_{0}$, we have a
morphism
\begin{align*}
&\hspace{1.1cm}\overline{F_{0}}:CokerL_{0}\rightarrow CokerM_{0}\\
&\hspace{3.1cm}X\mapsto F_{0}(X),\\
&\hspace{1.7cm}X\xrightarrow[]{[A,f]}Y\mapsto
F_0(X)\xrightarrow[]{[F(A),g]}F_0(Y)
\end{align*}
where $g$ is the composition
$F_{0}(X)\xrightarrow[]{F_{0}(f)}F_0(L_{0}(A)+Y)\backsimeq
F_{0}L_{0}(A)+F_{0}(Y)\xrightarrow[]{(\varepsilon_{0})_{A}+1}M_{0}F(A)+F_{0}(Y)$.
Moreover, there is a commutative diagram
\begin{center}
\scalebox{0.9}[0.85]{\includegraphics{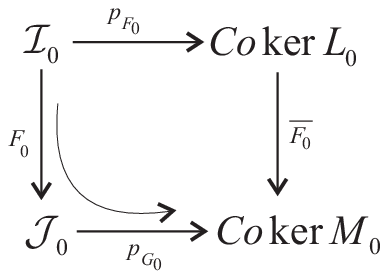}}
\end{center}
From the relative 2-exactness of injective resolution of
$\cR$-2-module,
there exist faithful morphisms
$\overline{L_{1}}:CokerF_{0}\rightarrow\cI_{1} $,
$\overline{M_{1}}:CokerG_{0}\rightarrow\cJ_{1}$ and 2-morphisms
$\varphi_{1}:\overline{L_1}\circ p_{L_0} \Rightarrow L_1$,
$\psi_{1}:\overline{M_{1}}\circ p_{M_0}\Rightarrow M_{1}$,
respectively. Then there exist 1-morphism $F_1:\cI_1\rightarrow
\cJ_1$ and 2-morphism $\overline{\varepsilon_{1}}:F_1\circ
\overline{L_1}\Rightarrow \overline{M_1}\circ\overline{F_0}$ from
the injectivity of $\cJ_{1}$.
From $\overline{\varepsilon_{1}}$ and $\overline{F_0}\circ p_{L_0}
=p_{M_0}\circ F_0$, we can define a 2-morphism
$\varepsilon_1:F_1\circ L_1\Rightarrow M_1\circ F_0$ by $ F_1\circ
L_1\xrightarrow[]{F_{1}(\varphi_{1}^{-1})}F_{1}\circ\overline{L_{1}}\circ
p_{L_0}\xrightarrow[]{\overline{\varepsilon_{1}}p_{L_0}}\overline{M_{1}}\circ\overline{F_0}\circ
p_{L_0}=\overline{M_1}\circ p_{M_0}\circ
F_0\xrightarrow[]{\psi_{1}F_0}M_1\circ F_0$, which is compatible
with $\varepsilon_{0}$.

Next we will construct $F_{n}$ and $\varepsilon_{n}:F_{n}\circ
L_{n}\Rightarrow M_{n}\circ F_{n-1}$ by induction on $n$.
Inductively, suppose $F_{i}$ and $\varepsilon_{i}$ have been
constructed for $i\leq n$ satisfying the compatible conditions.
Consider the morphism $F_{n}:\cI_{n}\rightarrow\cJ_{n}$, there is an
induced morphism
\begin{align*}
&\overline{F_{n}}:Coker(\alpha_{n},L_{n})\rightarrow
Coker(\beta_{n},M_{n})\\
&\hspace{2.9cm}X\mapsto F_{n}(X),\\
&\hspace{0.4cm}X\xrightarrow[]{[X_{n-1},x_{n-1}]}Y\mapsto
F_{n}(X)\xrightarrow[]{[F_{n-1}(X_{n-1}),y_{n-1}]}F_{n}(Y)
\end{align*}
where $y_{n-1}$ is the composition
$F_{n}(X)\xrightarrow[]{F_{n}(x_{n-1})}F_{n}(L_{n}(X_{n-1})+Y)\backsimeq
F_{n}L_{n}(X_{n-1})+F_{n}(Y)\xrightarrow[]{(\varepsilon_{n})_{X_{n-1}}+1}M_{n}F_{n-1}(X_{n-1})+F_{n}(Y)$.
Moreover, there is the following commutative diagram
\begin{center}
\scalebox{0.9}[0.85]{\includegraphics{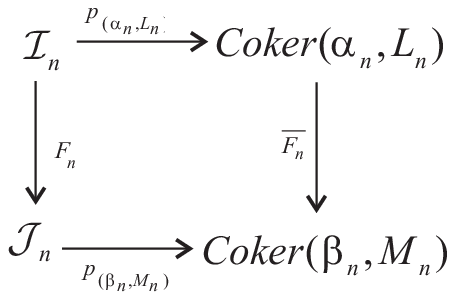}}
\end{center}

Using the relative 2-exactness of injective resolutions of $\cA$ and
$\cB$, we have the following diagram
\begin{center}
\scalebox{0.9}[0.85]{\includegraphics{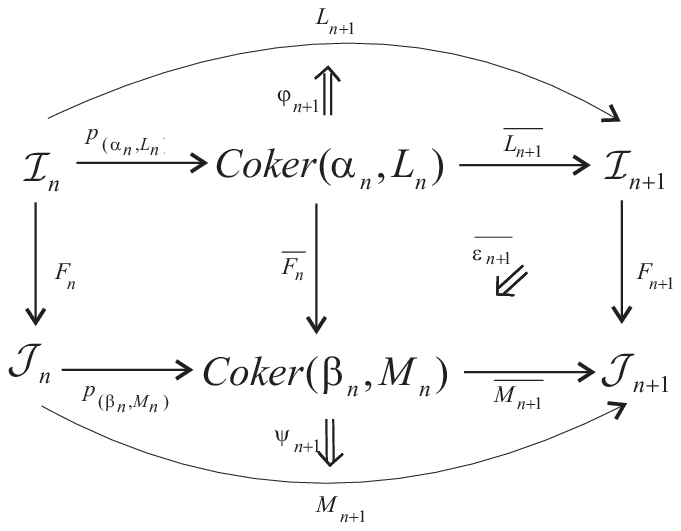}}
\end{center}
The existence of $F_{n+1}$ and $\overline{\varepsilon_{n+1}}$ come
from the injectivity of $\cJ_{n+1}$. Similar to the appearing of
$\varepsilon_{1}$, there is a 2-morphism $\varepsilon_{n+1}$ given
by $\overline{\varepsilon_{n+1}},\varphi_{n+1},\psi_{n+1}$, and is
compatible with $\varepsilon_{n}$.


Next, we show the uniqueness of $(F_{\cdot},\varepsilon_{\cdot})$ up
to 2-cochain homotopy in ($\cR$-2-Mod). Suppose
$(G_{\cdot},\mu_{\cdot})$ is another morphism of injective
resolutions. We will construct the 1-morphism
$H_{n}:\cI_{n+1}\rightarrow \cJ_{n},$ and 2-morphism $
\tau_{n}:F_n\Rightarrow M_{n}\circ H_{n-1}+H_{n}\circ L_{n+1}+G_n$
by induction on $n$. If $n<0$, $\cI_{n}=0$, so we get $H_n=0$. If
$n=0$, there is a 1-morphism $F_0-G_0:CokerL_0\rightarrow \cJ_0$(see
the following remark). From the faithful morphism
$\overline{L_1}:CokerL_0\rightarrow \cI_1$ and the injectivity of
$\cJ_0$, there exist a 1-morphism $H_0:\cI_1\rightarrow \cJ_0$ and
2-morphism $\tau_{0}^{'}:H_1\overline{L_1}\Rightarrow F_0-G_0$. Then
we get a 2-morphism $\tau_{0}:F_0\Rightarrow H_0\circ L_1+G_0$.

Inductively, we suppose given family of morphisms
$(H_i,\tau_i)_{i\leq n}$ so that
$H_i:\cI_i\rightarrow\cJ_{i-1},\tau_{i}:F_i\Rightarrow M_{i}\circ
H_{i-1}+H_{i}\circ L_{i+1}+G_{i}.$ Consider the 1-morphism
$F_n-G_n-M_{n}\circ H_{n-1}:Coker(\alpha_n,L_n)\rightarrow\cJ_n $
and the faithful morphism
$\overline{L_{n+1}}:Coker(\alpha_n,L_n)\rightarrow\cI_{n+1}$, the
injectivity of $\cJ_{n}$, there exist a 1-morphism
$H_n:\cI_{n+1}\rightarrow\cJ_{n}$ and a 2-morphism $\tau_{n}^{'}:
H_n\circ\overline{L_{n+1}}\Rightarrow F_n-G_n-M_{n}\circ H_{n-1}$.
Then we get a 2-morphism $\tau:F_n\Rightarrow M_{n}\circ
H_{n-1}+H_{n}\circ L_{n+1}+G_n$.
\end{proof}
\begin{Rek}
Morphism $F_0-G_0:CokerL_0\rightarrow \cJ_0$ is given by:

For an object $X$ in $CokerL_0$, define $(F_0-G_0)(X)$ by
$F_0(X)-G_0(X)$(see \cite{2,4}).

For a morphism $[A,f]:X\rightarrow Y$ in $CokerL_0$, where $A\in
obj(\cA),\ f:X\rightarrow L_0(A)+Y$, there are morphisms
$F_0(X)\xrightarrow[]{F_{0}(f)} F_0(L_0(A)+Y)\backsimeq
F_0L_0(A)+F_0(Y)\xrightarrow[]{(\varepsilon_{0})_{A}+1}M_0F(A)+F_0(Y)$,
$G_0(X)\xrightarrow[]{G_{0}(f)} G_0(L_0(A)+Y)\backsimeq
G_0L_0(A)+G_0(Y)\xrightarrow[]{(\mu_{0})_{A}+1}M_0F(A)+G_0(Y)$,
define $(F_0-G_0)([A,f]):F_0(X)-G_0(X)\rightarrow F_0(Y)-G_0(Y)$ by
the following morphisms under bifunctor $+$, i.e.
$(F_0-G_0)([A,f])\triangleq ((\varepsilon_{0})_{A}+1)\circ
F_0(f)-((\mu_{0})_{A}+1)\circ G_0(f)$.

$F_0-G_0$ is well-defined. In fact, if
$[A,f]=[A^{'},f^{'}]:X\rightarrow Y$ in $CokerL_0$, there exists a
morphism $a:A\rightarrow A^{'}$ such that the following diagram
commutes
\begin{center}
\scalebox{0.9}[0.85]{\includegraphics{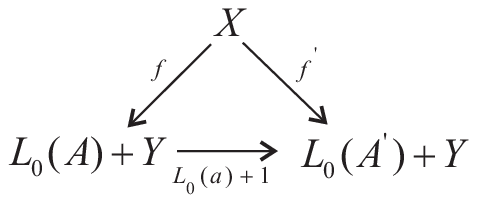}}
\end{center}
Give the operation of $F_0$ on the above commutative diagram,
together with the definitions of 1-morphism and 2-morphism in
($\cR$-2-Mod), we have the following commutative diagrams
\begin{center}
\scalebox{0.9}[0.85]{\includegraphics{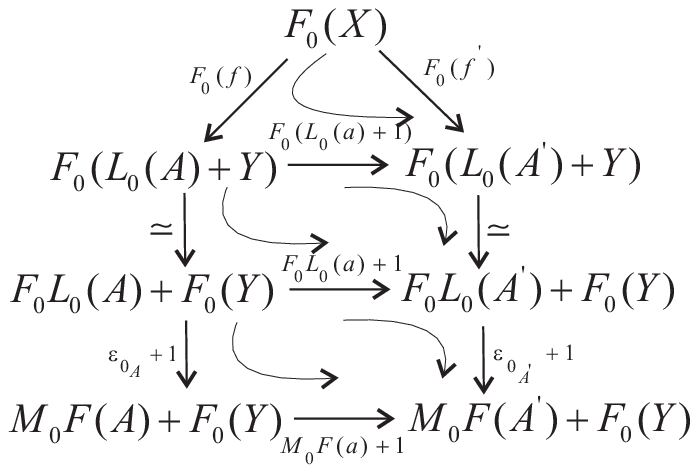}}
\end{center}
For $G_0$, there is a similar diagrams. Then we have
$(F_0-G_0)([A,f])=(F_0-G_0)([A^{'},f^{'}])$.

$F_0-G_0$ is a functor, i.e. for given morphisms
$X\xrightarrow[]{[A_1,f_1]}Y\xrightarrow[]{[A_2,f_2]}Z$ in
$CokerL_0$, there is
$(F_0-G_0)([A_2,f_2]\circ[A_1,f_1])=(F_0-G_0)([A_2,f_2])\circ(F_0-G_0)([A_1,f_1])$.
In fact, $[A_2,f_2]\circ[A_1,f_1]=[A_1+A_2,(1+f_2)\circ f_1]$,
$(F_0-G_0)([A_2,f_2]\circ[A_1,f_1])=(F_0-G_0)([A_1+A_2,(1+f_2)\circ
f_1])=((\varepsilon_0)_{A_1+A_2}+1)\circ F_0((1+f_2)\circ
f_1)-((\mu_0)_{A_1+A_2}+1)\circ G_0((1+f_2)\circ
f_1)=((\varepsilon_0)_{A_1}+(\varepsilon_0)_{A_2}+1)\circ
F_0(1+f_2)\circ F_0(f_1)-((\mu_0)_{A_1}+(\mu_0)_{A_2}+1)\circ
G_0(1+f_2)\circ G_0(f_1)=((\varepsilon_{0})_{A_2}+1)\circ
F_0(f_2)-((\mu_{0})_{A_2}+1)\circ
G_0(f_2)((\varepsilon_{0})_{A_1}+1)\circ
F_0(f_1)-((\mu_{0})_{A_1}+1)\circ
G_0(f_1)=(F_0-G_0)([A_2,f_2])\circ(F_0-G_0)([A_1,f_1]).$

$F_0-G_0$ is an $\cR$-homomorphism of $\cR$-2-modules from $F_0,\
G_0$ are.

Morphism $F_n-G_n-M_{n}H_{n-1}:Coker(\alpha_n,L_n)\rightarrow \cJ_n$
is given by the following way:

For an object $X$ in $Coker(\alpha_n,L_n)$, define
$(F_n-G_n-M_{n}H_{n-1})(X)$ by $F_n(X)-G_n(X)-M_{n}H_{n-1}(X)$(see
\cite{2,4}).

For a morphism $[X_{n-1},x_{n-1}]:X\rightarrow Y$ in
$Coker(\alpha_n,L_n)$, where $X_{n-1}\in obj(\cI_{n-1}),\
x_{n-1}:L_{n}(X_{n-1})+Y$. Consider the following composition
morphisms $F_n(X)\xrightarrow[]{F_{n}(x_{n-1})}
F_n(L_n(X_{n-1})+Y)\backsimeq
F_nL_n(X_{n-1})+F_n(Y)\xrightarrow[]{(\varepsilon_{n})_{X_{n-1}}+1}M_{n}F_{n-1}(X_{n-1})+F_n(Y)$,
$G_n(X)\xrightarrow[]{G_{n}(x_{n-1})}G_n(L_n(X_{n-1})+Y)\backsimeq
G_nL_n(X_{n-1})+G_n(Y)\xrightarrow[]{(\mu_{n})_{X_{n-1}}+1}M_{n}G_{n-1}(X_{n-1})+G_n(Y)$,
$M_nH_{n-1}(X)\xrightarrow[]{M_{n}H_{n-1}(x_{n-1})}
M_nH_{n-1}(L_n(X_{n-1})+Y)\backsimeq
M_nH_{n-1}L_n(X_{n-1})+M_nH_{n-1}(Y)$, together with 2-morphism
$\tau_{n-1}:F_{n-1}\Rightarrow M_{n-1}H_{n-2}+H_{n-1}L_{n}+G_{n-1}$,
define $(F_n-G_n-M_nH_{n-1})([X_{n-1},x_{n-1}])$ by the operation of
$+$ of above morphisms.

Similarly, $F_n-G_n-M_nH_{n-1}$ is an $\cR$-homomorphism of
$\cR$-2-modules.
\end{Rek}

\section{Derived 2-Functor in ($\cR$-2-Mod)}
In this section, we will define the right derived 2-functor between
the abelian 2-categories ($\cR$-2-Mod) and ($\cS$-2-Mod), which have
enough injective objects\cite{20,24}.

\begin{Def}Let $\cR,\cS$ be two 2-rings.
An additive 2-functor(\cite{2}) $T$:
($\cR$-2-Mod)$\rightarrow$($\cS$-2-Mod) is called left relative
2-exact if the relative 2-exactness of
\begin{center}
\scalebox{0.9}[0.85]{\includegraphics{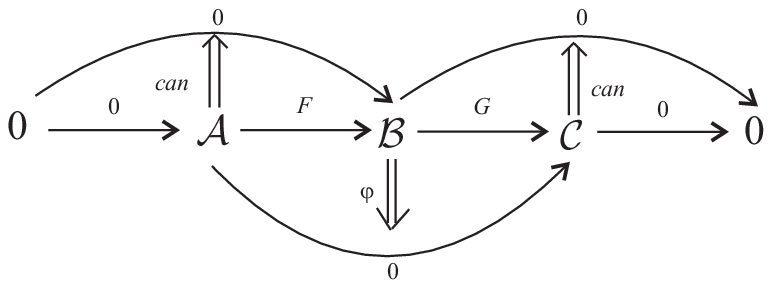}}
\end{center}
in $\cA,\cB$ and $\cC$ implies relative 2-exactness of
\begin{center}
\scalebox{0.9}[0.85]{\includegraphics{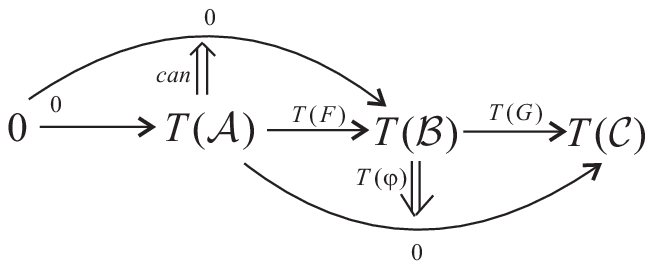}}
\end{center}
in $T(\cA)$ and $T(\cB)$.
\end{Def}

By Remark 2 and Proposition 1, Theorem 1, there is
\begin{cor}
Let $T$: ($\cR$-2-Mod)$\rightarrow$($\cS$-2-Mod) be an additive
2-functor,
and $\cA$ be an object of ($\cR$-2-Mod). For two injective
resolutions $\cI_{\cdot},\ \cJ_{\cdot}$ of $\cA$, there is an
equivalence between cohomology $\cR$-2-modules
$\cH^{\cdot}(T(\cI_{\cdot}))$ and $\cH^{\cdot}(T(\cJ_{\cdot}))$.
\end{cor}

Let $T$: ($\cR$-2-Mod)$\rightarrow$($\cS$-2-Mod) be an additive
2-functor. There is a 2-functor
\begin{align*}
&\cR^{i}T:(\cR\textrm{-2-Mod})\rightarrow (\cS\textrm{-2-Mod})\\
&\hspace{2.2cm}\cA\mapsto \cR^{i}T(\cA),\\
&\hspace{1.3cm}\cA\xrightarrow[]{F}\cB\mapsto
\cR^{i}T(\cA)\xrightarrow[]{\cR^{i}T(F)} \cR^{i}T(\cB)
\end{align*}
where $\cR^{i}T(\cA)$ is defined by $\cH^{i}(T(\cI_{\cdot}))$, and
$\cI_{\cdot}$ is the injective resolution of $\cA$. $\cR^{i}T$ is a
well-defined 2-functor from the properties of additive 2-functor and
Corollary 1.

\begin{cor}
Let $T$: ($\cR$-2-Mod))$\rightarrow$($\cR$-2-Mod)) be a left
relative 2-exact 2-functor, and $\cA$ be an injective object in
($\cR$-2-Mod)). Then $\cR^{i}T(\cA)=0$ for $i\neq0$.
\end{cor}

The following is a basic property of derived functors.
\begin{Thm}
The left derived 2-functor $\cL^{*}T$ takes
the sequence of $\cR$-2-modules
\begin{center}
\scalebox{0.9}[0.85]{\includegraphics{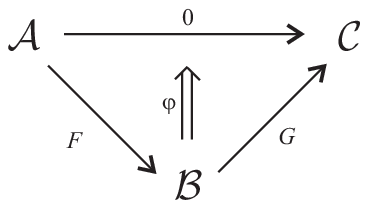}}
\end{center}
which is relative 2-exact in $\cA,\ \cB,\ \cC$ to a long sequence
2-exact(\cite{1,6})in each point in ($\cS$-2-Mod)
\begin{center}
\scalebox{0.9}[0.85]{\includegraphics{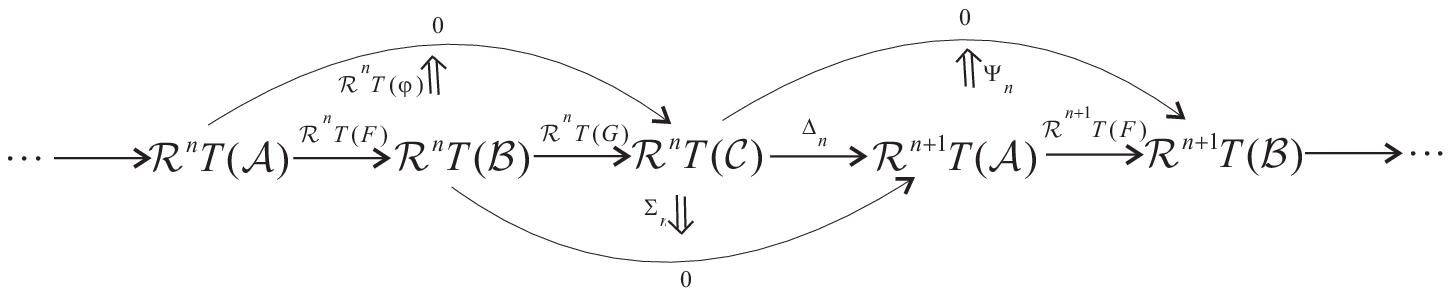}}
\end{center}
\end{Thm}

In order to prove this theorem, we need the following Lemmas.
\begin{Lem}
Let $\cI$ and $\cJ$ be injective objects in ($\cR$-2-Mod). Then the
product category $\cI\times\cJ$ is an injective object in
($\cR$-2-Mod).
\end{Lem}
\begin{proof}
First we know that $\cI\times\cJ$ is an $\cR$-2-module(\cite{4}). So
we need to prove the injectivity of it. There are canonical
morphisms
$$
\cI\xleftarrow[]{p_1}\cI\times\cJ\xrightarrow[]{p_2}\cJ,\
\cI\xrightarrow[]{i_1}\cI\times\cJ\xleftarrow[]{i_2}\cJ.
$$
For any morphism $G:\cA\rightarrow \cI\times\cJ$, there are
composition morphisms
$G_1:\cA\xrightarrow[]{G}\cI\times\cJ\xrightarrow[]{p_1}\cI$,
$G_2:\cB\xrightarrow[]{G}\cI\times\cJ\xrightarrow[]{p_2}\cJ$. Then
for a faithful morphism $F:\cA\rightarrow \cB$, there exist
1-morphisms $G_{1}^{'}:\cB\rightarrow \cI$,
$G_{2}^{'}:\cB\rightarrow \cJ$ and 2-morphisms $h_{1}:G_{1}^{'}\circ
F\Rightarrow G_{1}$, $h_{2}:G_{2}^{'}\circ F \Rightarrow G_{2}$
since $\cI$ and $\cJ$ are injective objects in ($\cR$-2-Mod).

So there are 1-morphism $G^{'}:\cB\rightarrow\cI\times\cJ$ given by
$G^{'}\triangleq (i_{1}G_{1}^{'},i_{2}G_{2}^{'})$ and 2-morphism
$h:G^{'}\circ F\Rightarrow G:\cA\rightarrow\cI\times\cJ$ given by
the composition
$h_{A}:G^{'}F(A)=(i_{1}G_{1}^{'},i_{2}G_{2}^{'})(F(A))=(i_{1}G_{1}^{'}F(A),i_{2}G_{2}^{'}F(A))\backsimeq
(G_{1}^{'}F(A),G_{2}^{'}F(A))\xrightarrow[]{((h_{1})_{A},(h_{2})_{A})}(G_{1}(A),G_{2}(A))=G(A)
$, for any $A\in obj(\cA)$.

Then $\cI\times\cJ$ is an injective object in ($\cR$-2-Mod).
\end{proof}

\begin{Lem}
Let $(F,\varphi,G):\cA\rightarrow\cB\rightarrow\cC$ be an extension
of $\cR$-2-modules in ($\cR$-2-Mod)(similar to the symmetric 2-group
case in \cite{11,2,6}), $(\cI_{\cdot},L_{\cdot},\alpha_{\cdot})$
$(\cJ_{\cdot},N_{\cdot},\beta_{\cdot})$ be injective resolutions of
$\cA$ and $\cC$, respectively. Then there is an injective resolution
$(\cK_{\cdot},M_{\cdot},\gamma_{\cdot})$ of $\cB$, such that
$\cI_{\cdot}\rightarrow \cK_{\cdot}\rightarrow\cJ_{\cdot}$ forms an
extension of 2-cochain complexes in ($\cR$-2-Mod).
\end{Lem}
\begin{proof}
We give the construction of injective resolution
$(\cK_{\cdot},M_{\cdot},\gamma_{\cdot})$ of $\cB$ in the following
steps.

Step 1. Since $\cI_0$ is an injective object, together with faithful
morphism $F:\cA\rightarrow\cB$ and 1-morphism $L_0:\cA\rightarrow
\cI_{0}$, there exist 1-morphism $\overline{L_0}:\cB\rightarrow
\cI_0$ and 2-morphism $h_0:\overline{L_0}\circ F\Rightarrow L_0$.
Then we can define a 1-morphism
\begin{align*}
&\hspace{2cm}M_0:\cB\rightarrow \cI_{0}\times\cJ_{0}\\
&\hspace{2.9cm}B\mapsto M_{0}(B)\triangleq
(\overline{L_{0}}(B),N_0G(B)),\\
&\hspace{1cm} (B_1)\xrightarrow[]{g}(B_2) \mapsto
(\overline{L_0}(g),N_0G(g)).
\end{align*}
Moreover, $M_0$ is faithful. In fact, if for any two morphisms
$g_1,g_2:B_1\rightarrow B_2$ in $obj(\cB)$, such that
$M_0(g_1)=M_0(g_2):M_0(B_1)\rightarrow M_0(B_2)$, i.e.
$(\overline{L_0}(g_1),N_0G(g_1))=(\overline{L_0}(g_2),N_0G(g_2))$.
From the definition of $\cI_0\times\cJ_0$, we have
$N_0G(g_1)=N_0G(g_2)$ and from $M_0$ is faithful, there is
$G(g_1)=G(g_2):G(B_1)\rightarrow G(B_2)$. By the universal property
of cokernel and the definition of extension of
$\cR$-2-modules(similar as symmetric 2-group case in \cite{11,2,6}),
there are a full and faithful morphism $G_0:CokerF\rightarrow\cC$
and a 2-morphism $\psi_0:G_0\circ p_{F}\Rightarrow
G:\cB\rightarrow\cC$. Then we have the identity
$G_0p_F(g_1)=G_0p_F(g_2):G(B_1)\rightarrow G(B_2)$ following from
the properties of $\psi_0$, i.e. $[0,g_1]=[0,g_2]:B_1\rightarrow
B_2$ in $CokerF$, so $g_1=g_2$.

Also,
\begin{center}
\scalebox{0.9}[0.85]{\includegraphics{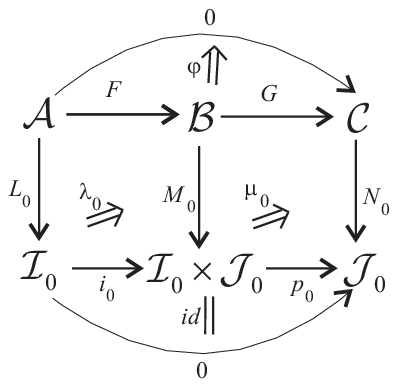}}
\end{center}
is the morphism of extensions in ($\cR$-2-Mod), where
$\lambda_0:i_0\circ L_0\Rightarrow M_{0}\circ F$ is given by
$h_{0}$, $\mu_{0}:p_0\circ M_0\Rightarrow N_0\circ G$ is the
identity.

Step 2. From the definition of relative 2-exactness, there are
faithful 1-morphisms $\overline{L_{1}}:CokerL_{0}\rightarrow\cI_{1}
$, $\overline{N_{1}}:CokerN_{0}\rightarrow \cJ_{1}$ as in the
following diagram
\begin{center}
\scalebox{0.9}[0.85]{\includegraphics{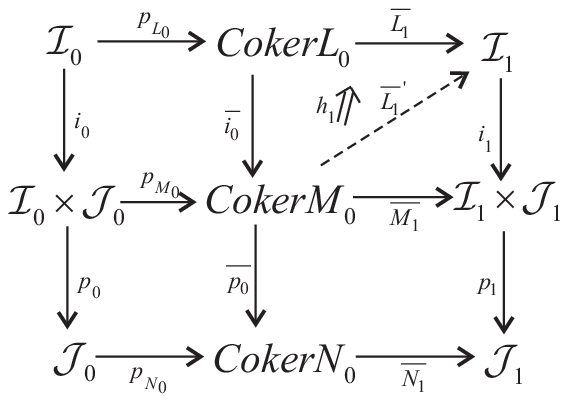}}
\end{center}
where $\overline{M_{1}}:CokerM_{0}\rightarrow\cI_{1}\times\cJ_{1}$
is given by $\overline{M_{1}}(X_0,Y_0)\triangleq
(\overline{L_{1}}^{'}(X_0,Y_0),
\overline{N_1}\overline{p_0}(X_0,Y_0)$, for any $(X_0,Y_0)\in
obj(\cI_{0}\times\cJ_{0})$, which is faithful from the proof of step
1.

Then we get a composition 1-morphism $M_{1}=\overline{M_{1}}\circ
p_{M_0}:\cI_{0}\times\cJ_{0}\rightarrow \cI_{1}\times\cJ_{1}$, and a
composition 2-morphism $\gamma_{1}:M_{1}\circ
M_{0}\Rightarrow\overline{M_1}p_{M_0}M_0\Rightarrow
\overline{M_{1}}\circ 0 \Rightarrow 0$, such that
\begin{center}
\scalebox{0.9}[0.85]{\includegraphics{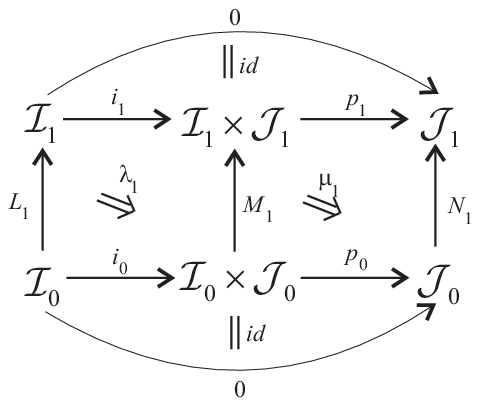}}
\end{center}
is a morphism of extensions in ($\cR$-2-Mod), where $\lambda_{1}$
and $\mu_{1}$ are given in the natural way as in step 1.

Step 3. From the definition of relative 2-exactness, there are
faithful 1-morphisms
$\overline{L_{2}}:Coker(\alpha_{1},L_{1})\rightarrow \cI_{2}$,
$\overline{N_{2}}:Coker(\beta_{1},N_{1})\rightarrow \cJ_{2}$ as in
the following diagram
\begin{center}
\scalebox{0.9}[0.85]{\includegraphics{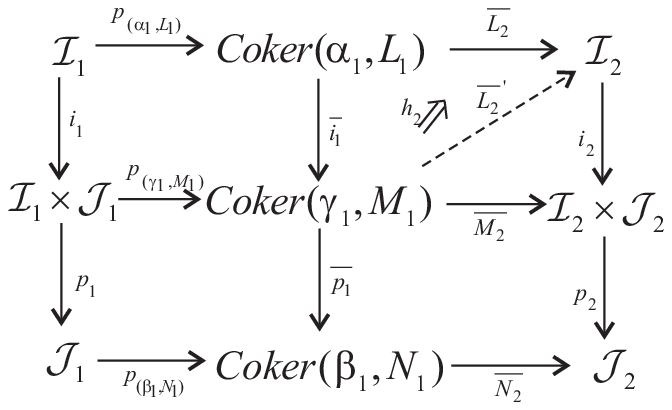}}
\end{center}
where $\overline{M_{2}}:Coker(\gamma_{1},M_{1})\rightarrow
\cI_{2}\times\cJ_{2}$ is given by the step 2, and also it is
faithful.


Then we get a composition 1-morphism $M_{2}=\overline{M_{2}}\circ
p_{(\gamma_{1},M_{1})}:\cI_{1}\times\cJ_{1}\rightarrow
\cI_{2}\times\cJ_{2}$, and a composition 2-morphism
$\gamma_{2}:M_{2}\circ M_{1}\Rightarrow \overline{M_{2}}\circ0
\Rightarrow 0$, such that
\begin{center}
\scalebox{0.9}[0.85]{\includegraphics{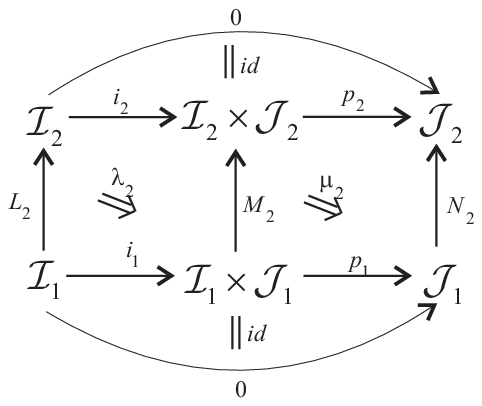}}
\end{center}
is a morphism of extensions in ($\cR$-2-Mod), where $\lambda_{2}$
and $\mu_{2}$ are given in the natural way as in step 1.

Using the same method, we get a complex
$(\cI\times\cJ_{\cdot},M_{\cdot},\gamma_{\cdot})$ of product
$\cR$-2-modules. Using the methods in Proposition 2, this complex is
relative 2-exact in each point, and
$(i_{\cdot},id,p_{\cdot}):\cI_{\cdot}\rightarrow
\cI_{\cdot}\times\cJ_{\cdot}\rightarrow\cJ_{\cdot}$ forms an
extension of complexes in ($\cR$-2-Mod).

Set $\cK_{n}=\cI_{n}\times\cJ_{n}$, for $n\geq 0$, which are
injective objects in ($\cR$-2-Mod) by Lemma 3.

\end{proof}

By the universal property of (bi)product of $\cR$-2-modules and the
property of additive 2-functor(\cite{2,4}). We get
\begin{Lem}
Let $T$: ($\cR$-2-Mod)$\rightarrow(\cS$-2-Mod) be an additive
2-functor, and $\cA,\ \cB$ be objects in ($\cR$-2-Mod). Then there
is an equivalence between $T(\cA\times\cB)$ and $T(\cA)\times
T(\cB)$ in ($\cS$-2-Mod).
\end{Lem}

Proof of Theorem 2. For $\cR$-2-modules $\cA$ and $\cC$, choose
injective resolutions $\cA\rightarrow\cI_{\cdot}$ and
$\cC\rightarrow\cJ_{\cdot}$. By Lemma 2 and Lemma 3, there is an
injective resolution $\cB\rightarrow \cI_{\cdot}\times\cJ_{\cdot}$
fitting into an extension
$\cI_{\cdot}\xrightarrow[]{i_{\cdot}}\cI_{\cdot}\times\cJ_{\cdot}\xrightarrow[]{p_{\cdot}}\cJ_{\cdot}$
of injective complexes in ($\cR$-2-Mod)(similar as \cite{1}). By
Lemma 4, we obtain a complexes of extension
$$T(\cI_{\cdot})\xrightarrow[]{T(i_{\cdot})}T(\cI_{\cdot}\times\cJ_{\cdot})\xrightarrow[]{T(p_{\cdot})}T(\cJ_{\cdot}).$$

Similar to the proof of Theorem 4.2 in \cite{11}, there is a long
sequence
\begin{center}
\scalebox{0.9}[0.85]{\includegraphics{p28.eps}}
\end{center}
which is 2-exact in each point.



\section*{Acknowledgements.} We would like to give our special thanks to Prof. Zhang-Ju LIU,
Prof. Yun-He SHENG for very helpful comments. We also thank Prof. Ke
WU and Prof. Shi-Kun WANG for useful discussions.

\bibliographystyle{model1b-num-names}

\noindent Fang Huang, Shao-Han Chen, Wei Chen\\
Department of Mathematics\\
 South China University of
 Technology\\
 Guangzhou 510641, P. R. China

\noindent Zhu-Jun Zheng\\
Department of Mathematics\\
 South China University of
Technology\\
 Guangzhou 510641, P. R. China \\
 and\\
Institute of Mathematics\\
Henan University\\  Kaifeng 475001, P. R.
China\\
E-mail: zhengzj@scut.edu.cn

\end{document}